\documentclass[twoside,11 pt]{article}%
\usepackage{amsmath}
\usepackage{graphicx}
\usepackage{url}
\usepackage{amsfonts}
\usepackage{amssymb}%
\usepackage{color}

\DeclareMathAlphabet{\cyrm}{U}{UWCyr}{m}{n}
\DeclareSymbolFont{cyrm}{U}{UWCyr}{m}{n}
\DeclareSymbolFontAlphabet{\cyrm}{cyrm}
\DeclareMathSymbol{\Evo}{\cyrm}{cyrm}{"03}
\newcommand{\Ker}{ \mathrm{Ker}}
\newcommand{\Rad}{ \mathrm{Rad}}

\newtheorem{theorem}{Theorem}[section]

\newtheorem{corollary}[theorem]{Corollary}

\newtheorem{definition}[theorem]{Definition}
\newtheorem{example}[theorem]{Example}

\newtheorem{lemma}[theorem]{Lemma}

\newtheorem{proposition}[theorem]{Proposition}
\newtheorem{remark}[theorem]{Remark}

\newenvironment{proof}[1][Proof]{\noindent\textbf{#1.} }{\
\rule{0.5em}{0.5em}}
\date{}

\begin{document}

\title{Contact geometry of multidimensional Monge-Amp\`{e}re equations: characteristics, intermediate integrals and  solutions}
\author{Dmitri Alekseevsky\thanks{School of Mathematics and Maxwell Institute for Mathematical Sciences, The Kings Buildings, JCMB, University of Edinburgh,
Mayfield  Road, Edinburgh, EH9 3JZ, UK,
D.Aleksee@ed.ac.uk}, Ricardo Alonso-Blanco\thanks{Departamento de Matem\'{a}ticas, Universidad de Salamanca, plaza de la Merced 1-4, 37008 Salamanca, Spain, ricardo@usal.es}, Gianni Manno\thanks{Dipartimento di Matematica e Applicazioni, Universit\`a degli Studi di Milano-Bicocca, via Cozzi 53, 20125 Milano, Italy, gianni.manno@unimib.it}, Fabrizio Pugliese\thanks{Dipartimento di
Matematica e Informatica,  Universit\`{a} di Salerno, via Ponte
don Melillo, 84084 Fisciano (SA), Italy, fpugliese@unisa.it}}
\maketitle

\begin{abstract}
\noindent
We study the geometry of multidimensional scalar $2^{nd}$ order PDEs (i.e. PDEs with $n$ independent variables) with one unknown function, viewed as hypersurfaces $\mathcal{E}$ in the Lagrangian Grassmann bundle $M^{(1)}$ over a $(2n+1)$-dimensional contact manifold $(M,\mathcal{C})$. We develop the theory of characteristics of the equation $\mathcal{E}$ in terms of contact geometry and of the geometry of Lagrangian Grassmannian and study their relationship with intermediate integrals of $\mathcal{E}$. After specifying the results to general Monge-Amp\`ere equations (MAEs), we focus our attention to MAEs of type introduced by Goursat in \cite{Goursat2}, i.e. MAEs of the form
\[
\det\left\|\frac{\partial^2 f}{\partial x^i\partial x^j}-b_{ij}\left(x,f,\nabla f\right)\right\|=0.
\]
We show that any MAE of the aforementioned class is associated with an
$n$-dimensional subdistribution $\mathcal{D}$ of the contact
distribution $\mathcal{C}$, and viceversa. We characterize this Goursat-type equations together with its intermediate integrals in terms of their characteristics and give a criterion of local contact equivalence. Finally, we develop a method of solutions of a Cauchy problem,
provided the existence of a suitable number of intermediate
integrals.
\end{abstract}

\medskip
\noindent\textbf{MSC Classification 2010:} 53D10, 35A30, 58A30, 58A17


\tableofcontents


\section{Introduction}

\subsection{Characteristics of PDEs, Cauchy-Kowalewski theorem and MAEs}

Characteristics of PDEs are a classic subject
(\cite{Goursat,Goursat2,Petro,Valiron}) as they are related to the
local existence and uniqueness of solutions of Cauchy problems. As an example, if
\begin{equation}\label{ equation}
F(x^{1},\dots,x^{n},z,p_{1},\dots,p_{n},p_{11},p_{12},\dots
p_{nn})=0
\end{equation}
where $z=z(x^1,\dots,x^n)$,
$p_i = \partial z/ \partial x^i, \,\, p_{ij} =
\partial^2 z/ \partial x^i \partial x^j$ is a scalar second order partial differential equation ($2^{nd}$ order PDE), the Cauchy problem consists of finding a solution  $z=f(x^1,\dots,x^n)$ of \eqref{ equation} which satisfies the following
conditions
\begin{equation}\label{Cauchy_problem}
f|_{(X^1(\mathbf{t}),\dots,X^n(\mathbf{t}))}=Z(\mathbf{t})\,,\,\,\,\left.
\frac{\partial f}{\partial x^{i}}\right|
_{(X^1(\mathbf{t}),\dots,X^n(\mathbf{t}))}=P_{i}(\mathbf{t}),
\end{equation}
where
\begin{equation}\label{Cauchy_datum}
\Phi(\mathbf{t})=(X^{1}(\mathbf{t}),\dots,X^n(\mathbf{t}),Z(\mathbf{t}),P_{1}(\mathbf{t}),\dots,P_n(\mathbf{t}))\,,
\,\,\, \mathbf{t}=(t_1,\dots,t_{n-1})
\end{equation}
is a given $(n-1)$-dimensional manifold, i.e. a \emph{Cauchy datum}; obviously, in (3) the choice of the parametrization is irrelevant. If Cauchy datum \eqref{Cauchy_datum} is \emph{non-characteristic}, then, in the $C^\infty$ case, Cauchy problem \eqref{Cauchy_problem} for Equation \eqref{ equation} admits, locally, a unique formal solution: in fact in this hypothesis we can put Equation \eqref{ equation} in the Cauchy-Kowalewski form (see Section \ref{sec.formal.solutions} for a geometric description). Under the same hypothesis, in the analytic case it admits a locally unique solution.

\smallskip\noindent
In the case $n=2$, non-characteristicity condition means that tangent
direction $v=\dot{\Phi }(0)$ at a point
$m=\Phi(0)=(\overline{x}^1,\overline{x}^2,\overline{z},\overline{p}_1,\overline{p}_2)$ of the ($1$-dimensional, in this case) Cauchy datum
satisfies the condition
\begin{equation}\label{non_characteristic}
\left.\frac{\partial F}{\partial p_{11}}\right|_{m^1}{v^2}^2 - \left.\frac{\partial F}{\partial p_{12}}\right|_{m^1}{v^1v^2} +
\left.\frac{\partial F}{\partial p_{22}}\right|_{m^1}{v^1}^2\neq 0
\end{equation}
for each
$m^1=(\overline{x}^1,\overline{x}^2,\overline{z},\overline{p}_1,\overline{p}_2,
\overline{p}_{11},\overline{p}_{12},\overline{p}_{22})$ satisfying \eqref{ equation}, where
\begin{equation*}
v=v^1({\partial}_{x^1} + \overline{p}_1\partial_z +
\overline{p}_{1i}\partial_{p_i}) + v^2({\partial}_{x^2} +
\overline{p}_2\partial_z + \overline{p}_{2i}\partial_{p_i}).
\end{equation*}
The  vector $v$ can be considered as an ``infinitesimal
Cauchy datum''.

\smallskip\noindent
From Equation \eqref{non_characteristic} it is clear that one can associate with any point $m^1$ satisfying \eqref{ equation} two (possibly imaginary) directions in the space $(x^i,z,p_i)$, namely, those annihilating \eqref{non_characteristic} (``characteristic lines''); if we let this point vary keeping the point $m$ fixed, these two directions form, in general, two cones at $m$. It is proved that the only PDEs for which these two cones degenerates in two $2$-dimensional planes are classical Monge-Amp\`ere equations (MAEs) (see for instance \cite{AMP_Acta,AMP_DGA}).

\smallskip\noindent
One of the targets of this paper is to see if a similar phenomenon occurs also in the case of MAEs with an arbitrary number of independent variables, which, of course, is considerably more complicated.

\smallskip\noindent
In fact, MAEs for $n=2$ have been intensely studied since the second half of XIX century by many \emph{g\'eom\`etres}, among them Darboux, Lie, Goursat (a systematic account of such investigations can be found in \cite{For} and \cite{Goursat}); later, this classical approach was put aside in favour of more ``hard analysis'' techniques. The last 40 years have witnessed a renewed interest in the differential-geometric approach to MAE's, mainly due to Lychagin and his school (see \cite{K2} and \cite{KLR} for an exhaustive bibliography). However, such results are focused on the classical case ($n=2$). Up to now, no serious effort has been made to extend the classical theory to the general multidimensional case (only very special cases have been studied). In fact, the main achievements so far obtained in this direction are due to Boillat and Lychagin.

\smallskip\noindent
Boillat \cite{Boillat1} noticed that MAEs with two independent variables were the only second order PDEs which are exceptional in the sense of Lax \cite{Lax}. This physical property was used in \cite{Ruggeri} to find the general form of a MAE in three independent variables, and in \cite{Boillat2} for the case of arbitrary independent variables. The result is that such general form is
\begin{equation}\label{eq.lin.comb.minors}
M_n+M_{n-1}+\dots M_0=0
\end{equation}
where $M_k$ is a linear combination (with functions of $x^i,z,p_i$ as coefficients) of all $k\times k$ minors of the Hessian matrix $\|z_{x^ix^j}\|$.

\smallskip\noindent
In \cite{Lychagin contact}, by introducing a new approach based on contact geometry, Lychagin defined multidimensional MAEs as the zero locus of a differential operator associated with a class of $n$-differential forms on a contact manifold. Locally, such PDEs are described by \eqref{eq.lin.comb.minors}. In the rest of the paper, when we write ``general MAEs'' we mean ``multidimensional MAEs in the sense of Lychagin''.

\smallskip\noindent
The oldest paper regarding the multidimensional generalization of the concept of MAEs dates back to Goursat. In \cite{Goursat2} he noticed that classical MAEs ($n=2$) can be obtained by substituting $dp_1=p_{11}dx^1+p_{12}dx^2$ and $dp_2=p_{12}dx^1+p_{22}dx^2$ in the following system
$$
\left\{
\begin{array}{l}
dp_1-b_{11}dx^1-b_{12}dx^2=0\\
dp_2-b_{21}dx^1-b_{22}dx^2=0 \qquad b_{ij}=b_{ij}(x^1,x^2,z,p_1,p_2)
\end{array}
\right.
$$
and by requiring its (non trivial) compatibility. Obviously, such ``horizontalization'' of the above Pfaffian system can be extended to any number $n$ of independent variables; namely, one can consider the system
$$
dp_i-\sum_{j=1}^n b_{ij}dx^j=0\,,\,\,\,i=1,\dots,n \,,\,\,\, b_{ij}=b_{ij}(x^1,\dots,x^n,z,p_1,\dots,p_n)\,,
$$
``horizontalize'' it $(dp_i=p_{ij}dx^j)$ and impose the compatibility condition, thus getting MAE
\begin{equation}\label{eq.MAEs.Goursat.type}
\det||p_{ij}-b_{ij}||=0.
\end{equation}
It turns out that the class of PDEs considered by Goursat is a subclass of those considered by Lychagin.

\smallskip\noindent
The above analytical procedure has a natural geometrical meaning, tightly linked with the fundamental notion of characteristics of a PDE. Such a connection, which was already studied in \cite{AMP_Acta,AMP_DGA} for $n=2$, will be extended below to the case of any number of independent variables. As we shall see, for $n>2$ the complexity of the problem drastically increases. For this purpose, as a first step we develop a coordinate free setting to the theory of characteristics of scalar second order PDEs (with $n$ independent variables) in terms of contact manifolds and Lagrangian Grassmannians, which we summarize below.

\smallskip\noindent
Let $(M,\mathcal{C})$ be a \emph{$(2n+1)$-dimensional contact
manifold}, i.e. a $(2n+1)$-dimensional manifold where $\mathcal{C}$ is a completely non integrable
distribution of codimension $1$. Locally $\mathcal{C}$ is the kernel of (a contact) $1$-form $\theta$ (which is defined up to a conformal factor) which in appropriate (contact or Darboux) coordinates $(x^i,z,p_i), i= 1,\dots, n$ has the
form
$$
\theta = dz - p_i dx^i.
$$
The restriction
$$
\omega = d \theta|_{\mathcal{C}}
$$
defines on each hyperplane $\mathcal{C}_m$ a conformal symplectic structure, of fundamental importance in contact geometry: in fact, Lagrangian (i.e. maximally $\omega$-isotropic) planes of $\mathcal{C}_m$ are tangent to maximal integral submanifolds of $\mathcal{C}$ and thus $n$-dimensional; for this reason, such submanifolds of $M$ are called \emph{Lagrangian} (or also \emph{Legendrian}). We denote by $\mathcal{L}(\mathcal{C}_m)$ the \emph{Grassmannian of Lagrangian planes} of $\mathcal{C}_m$  and by
$$
\pi: M^{(1)} = \bigcup_{m \in M}\mathcal{L}(\mathcal{C}_m) \to M $$
the  bundle of Lagrangian planes. Contact coordinates $(x^i,z,p_i)$ on $M$ induce coordinates on $M^{(1)}$: a point $m^1 \equiv L_{m^1} \in M^{(1)}$
has coordinates $(x^i, z, p_i, p_{ij})$, $1 \leq i \leq j \leq n$ iff the
corresponding Lagrangian plane $L_{m^1}$ is given by:
\begin{equation*}
m^1\equiv L_{m^1} = \langle\widehat{\partial}_{x^i} + p_{ij}\partial_{p_j}\rangle\,,\,\, \widehat{\partial}_{x^i}\overset{\textrm{def}}{=} \partial_{x^i}+ p_i \partial_z
\end{equation*}
with $\|p_{ij}\|$ a symmetric matrix.

\smallskip\noindent
A \emph{scalar $2^{nd}$ order PDE with $n$ independent variables with one unknown function} is defined as a
hypersurface $\mathcal{E}$ of $M^{(1)}$ and its \emph{solutions} are Lagrangian submanifolds $\Sigma \subset M$ such that $T\Sigma \subset \mathcal{E}$. In view of reasonings made at the beginning of the section, a \emph{Cauchy datum} for $\mathcal{E}$ is defined simply as an $(n-1)$-dimensional submanifold of $M$ which in view of \eqref{Cauchy_problem} must be also integral of $\mathcal{C}$. The restriction on $\mathcal{E}$ of fibre bundle $\pi$ is a bundle over $M$ whose fibre at $m$ is denoted by $\mathcal{E}_m$:
\begin{equation}\label{eq.principal}
\mathcal{E}_m:=\mathcal{E}\cap \mathcal{L}(\mathcal{C}_m).
\end{equation}
$\mathcal{E}_m$ is a hypersurface of the Grassmannian  $\mathcal{L}(\mathcal{C}_m)$ of Lagrangian planes of $\mathcal{C}$. A straightforward computation shows that the set of Lagrangian planes at $m\in M$ containing a given $(n-1)$-dimensional isotropic subspace is a curve in $\mathcal{L}(\mathcal{C}_m)$: condition \eqref{non_characteristic} (in the case $n=2$) means that the curve formed by Lagrangian planes containing $v$ is not tangent to $\mathcal{E}_m$ at $m^1$. This condition can be easily generalized to any dimension: we can define a \emph{characteristic subspace for $\mathcal{E}$ at $m^1$} as a hyperplane of $L_{m^1}$ such that the curve in $\mathcal{L}(\mathcal{C}_m)$ whose points are Lagrangian planes containing it is tangent to $\mathcal{E}_m$ at $m^1$. The tangent space to this curve at $m^1$ is called a \emph{characteristic direction for $\mathcal{E}$ at $m^1$}.

\smallskip\noindent
By means of previous geometric concepts, we are able to give an intrinsic definition of MAEs of form \eqref{eq.lin.comb.minors} and \eqref{eq.MAEs.Goursat.type}.
The former describe, locally, hypersurfaces $\mathcal{E}_\Omega$ of $M^{(1)}$ formed by Lagrangian planes which annihilate an $n$-form $\Omega$ on $M$:
\begin{equation}\label{eq.MAE.intr.Omega}
\mathcal{E}_{\Omega}=\{{m}^1\in M^{(1)}
\,\,\big|\,\,\Omega|_{L_{{m}^1}}=0\},
\end{equation}
whereas the latter hypersurfaces $\mathcal{E}_{\mathcal{D}}$ of $M^{(1)}$ formed by Lagrangian planes
which non trivially intersect an $n$-dimensional subdistribution $\mathcal{D}$ of $\mathcal{C}$:
\begin{equation}\label{eq.MAE.intr.D}
\mathcal{E}_{\mathcal{D}}=\{ m^1\in M^{(1)} \,\,\big|\,\, L_{m^1}\cap \mathcal{D}_{\pi(m^1)}\neq 0  \}.
\end{equation}
It is easy to realize that MAEs of type $\mathcal{E}_{\mathcal{D}}$ are associated with decomposable $n$-forms on $M$.

\subsection{Main results and description of the paper}

All we said so far shows that characteristics of a PDE $\mathcal{E}$ are of ``point'' nature, in the sense that any information regarding them is contained in their fibres \eqref{eq.principal}. This justifies the importance of studying conformal properties of the Grassmannian of Lagrangian planes $\mathcal{L}(V)$ of a generic symplectic space $(V,\omega)$ together with its submanifolds. In \cite{FHK} an interpretation of special MAEs with constant coefficients is given in terms of Lagrangian Grassmannians. We concentrate mostly on hypersurfaces of Lagrangian Grassmannians, as the fibre \eqref{eq.principal} of a PDE is a hypersurface of $\mathcal{L}(\mathcal{C}_m)$. We study these subjects in Sections \ref{sec.Lagr.Grass.of.sympl.sp} and \ref{sec.sub.of.Lagr.Grass}, and then we reformulate the results in the languages of PDEs and MAEs in Section \ref{subsec.PDE}.

\smallskip\noindent
In Section \ref{sec.Lagr.Grass.of.sympl.sp} we describe the main geometric structures of the Lagrangian Grassmannian $\mathcal{L}(V)$.
We denote by $\mathcal{T}(\mathcal{L}(V))$ the \emph{tautological vector bundle of $\mathcal{L}(V)$}, i.e. the vector bundle on $\mathcal{L}(V)$ whose fibre at a point $L\in\mathcal{L}(V)$ is the vector space $L$. The main geometric structure of $\mathcal{L}(V)$ is  the ``symmetric Grassmann structure" i.e.  a canonical identification
\begin{equation}\label{eq.ident}
g:T\mathcal{L}(V) \overset{\sim}{\rightarrow} S^2 \big(\mathcal{T}^* (\mathcal{L}(V))\big)\,,\,\,\,v\mapsto g^v
\end{equation}
of the tangent bundle with the  symmetric square of the  dual tautological bundle. To keep the notation simple, we continue to denote the inverse of the dual map of \eqref{eq.ident} by $g$:
\begin{equation}\label{eq.ident.dual}
g:T^\ast\mathcal{L}(V) \overset{\sim}{\rightarrow} S^2 \big(\mathcal{T} (\mathcal{L}(V))\big)\,,\,\,\,\rho\mapsto g_\rho
\end{equation}
Note that there is no ambiguity in denoting by $g$ both the maps \eqref{eq.ident} and \eqref{eq.ident.dual} since vectors appear as superscripts whereas covectors as subscripts. Thus one can define the rank of vectors (resp. covectors) as the rank of the corresponding bilinear form through \eqref{eq.ident} (resp. \eqref{eq.ident.dual}). We underline that both $g^v$ and $g_\rho$ change conformally if the symplectic form $\omega$ change conformally.\\
The manifold $\mathcal{L}(V)$ has a natural Pl\"ucker embedding
into the projective space $\mathbb{P} \Lambda^n V$ so that any tangent vector
$\dot L \in T_L \mathcal{L}(V)$ defines a projective line
$\ell(L,\dot L) \subset \mathbb{P} \Lambda^n V$, that we show it belongs to
$\mathcal{L}(V)$ iff  $\textrm{rank}(\dot{L})=1$.

\smallskip\noindent
In Section \ref{sec.sub.of.Lagr.Grass} we study geometry of  submanifolds (mostly, hypersurfaces) of
$\mathcal{L}(V)$.\\
In view of \eqref{eq.ident.dual}, with any hypersurface $\textrm{E} =\{F(p_{ij}) =0\}$ of $\mathcal{L}(V)$ it is associated the (possibly degenerate) conformal metric
$$
g_{\textrm{E}}=[g_{dF}|_{\mathrm{E}}]\,,
$$
which turns out to be independent of the function $F$. Characteristic subspaces and characteristic directions of \textrm{E} are defined as follows.
Any subspace $U \subset V$ defines a distinguished submanifold $U^{(1)}$ of $\mathcal{L}(V)$, which we call the \emph{(first) prolongation of $U$}, formed by Lagrangian planes containing $U$ if $\dim U \leq n$ or which are contained in $U$ otherwise. An isotropic subspace $U\subset L\in\mathrm{E}$ is called a \emph{characteristic subspace} for
$\textrm{E}$ at $L$ if $U^{(1)}$ is tangent to $\textrm{E}$ at $L$. In the case that $U$ is an $(n-1)$-dimensional characteristic subspace for \textrm{E} at $L$, $U^{(1)}$ is $1$-dimensional and the tangent space $T_LU^{(1)}$ is called a \emph{characteristic direction} (for \textrm{E} at $L$): its elements are vectors of rank $1$.\\
The converse is also true: the radical of $g^{\dot{L}}$ (see \eqref{eq.ident}) where $\dot{L}$ spans a characteristic direction for $\textrm{E}$ at $L$ (i.e. $\dot{L}$ is a vector of $T_L\textrm{E}$ of rank $1$) is a characteristic subspace for
$\textrm{E}$ at $L\in\textrm{E}$. In other words, the projective line $\ell(L,\dot{L})$ associated with such $\dot{L}$ is tangent to $\textrm{E}$ (via the Pl\"ucker embedding). Up to sign, $g^{\dot L} = \eta \otimes \eta$ where $\eta \in
L^*$ is a $g_{\textrm{E}}$-isotropic covector.

\noindent
An important class of hypersurfaces  of $\mathcal{L}(V) \subset \mathbb{P}\Lambda^nV$ is that of \emph{hyperplane sections} of $\mathbb{P}\Lambda^n(V)$: they are the intersection of $\mathcal{L}(V)$ with a hyperplane of $\mathbb{P}\Lambda^nV$ (via the Pl\"ucker embedding). Since any hyperplane of $\mathbb{P}\Lambda^nV$ is given by $\{\Omega=0\}$ where $\Omega \in \Lambda^nV^*$, we denote such a hypersurface by $\textrm{E}_\Omega$. Hypersurfaces of type $\textrm{E}_\Omega$ are the prototype of fibres \eqref{eq.principal} of a general MAE, i.e. of type \eqref{eq.MAE.intr.Omega}.

\noindent
At the end of this Section \ref{sec.sub.of.Lagr.Grass}, we study hypersurfaces $\textrm{E}_D$ associated with an $n$-plane $D \subset V$. By definition, such a hypersurface consists of Lagrangian planes which have non-trivial intersection with $D$.
It is easy to realize that these hypersurfaces are special hyperplane sections of $\mathbb{P}\Lambda^n(V)$: they are defined by decomposable $n$-forms on $\mathcal{L}(V)$. Hypersurfaces of type $\textrm{E}_D$ are the prototype of fibres \eqref{eq.principal} of a MAE of Goursat type, i.e. of type \eqref{eq.MAE.intr.D}.
\\
The main results of Section \ref{sec.sub.of.Lagr.Grass} can be summarized as follows:
\begin{itemize}
\item Characteristic subspaces for a hypersurface $\textrm{E}$ of $\mathcal{L}(V)$ are those whose annihilator is $g_{\textrm{E}}$-isotropic (Theorem \ref{th.fab}). By using this, we find a relationship between the decomposability of $g_{\textrm{E}}$ and the behavior of characteristic subspaces (Theorem \ref{th.rotating.debole});
\item The projective line $\ell(L,\dot{L})$ associated with a characteristic vector $\dot{L}$ of a hyperplane section $\textrm{E}_\Omega$ is included in $\textrm{E}_\Omega$ (we say that $\dot{L}$ is \emph{strongly characteristic}). In other word, if a hyperplane $H$ of $L\in\mathcal{L}(V)$ is characteristic at $L$ for a hypersurface of type $\textrm{E}_\Omega$, then it is characteristic for any $\overline{L}\in\textrm{E}_\Omega$ such that $\overline{L}\supset H$ (Theorem \ref{th.strongly.char}).  We also describe $H$ in terms of isotropy of $\Omega$ (Theorem \ref{th.car.hyp});
\item A hypersurface of type $\textrm{E}_D$ can be associated only with two $n$-dimensional planes of $V$ which are mutually symplectically orthogonal (Theorem \ref{th.D.equal.D.orth});
\item Conformal metric $g_{\textrm{E}_D}$ is decomposable: it has rank equal to $1$ if $D$ is Lagrangian and rank $2$ otherwise. For each regular point $L\in\textrm{E}_D$ we have that $(g_{\textrm{E}_D})_{L}=\ell_{L}\vee\ell'_{L}$, where $\ell_{L}=L\cap D$ and $\ell_{L}'=L\cap D^\perp$ are lines. Then we have the following correspondence:
    \[
    L\in\textrm{E}_D\longmapsto (\ell_{L},\ell'_{L}).
    \]
$\textrm{E}_D$ possesses two
$(n-2)$-parametric families $H$ and $H'$ of characteristic hyperplanes of $L$ which rotate, respectively, around the line $\ell_{L}$ and resp. $\ell_{L}'$: if we let vary the point $L$ on ${\textrm{E}_D}$, the corresponding lines fill the $n$-dimensional
space $D$ (resp. $D^\perp$). In other words, we can reconstruct $\textrm{E}_D$ starting from its characteristics (Theorem \ref{theorem.main.2}).
\end{itemize}
By substituting $\mathcal{L}(V)\leftrightarrow\mathcal{L}(\mathcal{C}_m)$, $\textrm{E}\leftrightarrow\mathcal{E}_m$, $\textrm{E}_\Omega\leftrightarrow (\mathcal{E}_\Omega)_m$, $\textrm{E}_D\leftrightarrow (\mathcal{E}_{\mathcal{D}})_m$, $g_{\textrm{E}}\leftrightarrow g_{(\mathcal{E}_m)}$ in the above points, we reformulate previous results in the language of PDEs in Sections \ref{subsec.PDE.parte.prima}, \ref{sec.char} and in that of MAEs in Sections \ref{sec.char.MAE}, \ref{sec.MAE.ass.with.D}.

\smallskip\noindent
In Section \ref{sec.ContactManifolds} we recall the basic notions of contact geometry and geometric theory of first order PDE. We also shortly describe the solution of the Cauchy problem by the method of characteristics.

\smallskip\noindent
In Section \ref{subsec.PDE}, beside the results that we described above, we give a criterion of local equivalence for a PDE to be a MAE of Goursat type (Theorem \eqref{th.criterion}).

\smallskip\noindent
For the sake of completeness, in Section \ref{sec.full.prolong} we deal with the full (or infinite) prolongation of a $2^{nd}$ order PDE.
We show that any $2^{nd}$ order PDE $\mathcal{E}$ is formally integrable provided that conformal metric $g_{\mathcal{E}}$ does not
vanish, and that a non-characteristic Cauchy problem has unique formal solution. In fact, finding necessary and sufficient conditions for the existence and uniqueness of the solution of the Cauchy problem is the historical motivation of the notion of characteristics.

\smallskip\noindent
In Section \ref{sec.int.int.PDEs.MAEs} we consider  intermediate integrals of $2^{nd}$ order PDEs with special attention to MAEs of type
$\mathcal{E}_{\mathcal{D}}$. The main results of the section are summarized below.
\begin{itemize}
\item The existence of an intermediate integral of a $2^{nd}$ order PDE is equivalent to the existence of a special vector field (Hamiltonian vector field) whose directions are strongly characteristic (Theorem \ref{th.f.int.int.Yf.str.char});

\item Intermediate integrals of $\mathcal{E}_{\mathcal{D}}$ coincide with the first integrals of the distribution $\mathcal{D}$ or $\mathcal{D^\perp}$ (Theorem \ref{th.equiv.int.int.with.first.int}). In particular, the existence of such a first integral implies the existence of a $C^\infty$ solution of $\mathcal{E}_{\mathcal{D}}$.

\item If $\mathcal{D}$ (or $\mathcal{D}^\perp$) possesses $n$ independent first integrals, we describe a method (going back to Monge and reinterpreted in contact geometric terms by Morimoto \cite{Morimoto}) of solution of any Cauchy problem associated with $\mathcal{E}_{\mathcal{D}}$ which involves only solutions of ordinary differential equations and finite equations (Theorem \ref{th.Monge.Morimoto}). We also show that, in this case, $\mathcal{E}_{\mathcal{D}}$ can be reconstructed by means of its intermediate integrals (Theorem \ref{th.ultimo}).
\end{itemize}

\medskip\noindent\textbf{\large Notations and conventions:}

\smallskip\noindent
In the rest of the paper Latin indices will run from $1$ to $n$, unless otherwise specified.
We will use Einstein convention.
We denote by $X\cdot\varrho$ the Lie derivative of a form $\varrho$ along a vector field $X$.
The symmetric tensor product will be denoted by $\vee$, i.e. $A\vee B=\frac{1}{2}(A\otimes B + B\otimes A)$.
The annihilator of a vector subspace $U$ will be denoted by $U^0$. We denote by $\langle v_i\rangle$ the linear span of vectors $v_1,\dots,v_n$.

\section{Geometry of  the Lagrangian Grassmannian $\mathcal{L}(V)$ }\label{sec.Lagr.Grass.of.sympl.sp}

\subsection{Lagrangian Grassmannian  $\mathcal{L}(V)$ and  its tautological bundle
$\mathcal{T}(\mathcal{L}(V))$}\label{sec.Lagr.Grass}

Let $(V,\omega)$ be a symplectic $2n$-dimensional vector space.
Recall that a  \textbf{Lagrangian plane}  is an isotropic
subspace $L \subset V$ of maximal dimension, i.e. an
$n$-dimensional subspace $L$ such that $\omega|_L =0$. We shall
denote by
$$
\mathcal{L}(V):=LGr(V)
$$
the \textbf{Grassmannian of Lagrangian planes} in $V$.

\smallskip\noindent
A smooth structure  of  the manifold  $\mathcal{L}(V)$   is
defined as follows. For any $L_{0}\in\mathcal{L}(V)$, we  choose a
complementary Lagrangian plane $L_{0}'\in\mathcal{L}(V)$, and a
symplectic basis $\{e_{i},e^i\}$ (i.e. $\omega(e_i,e^j)=\delta_i^j$) such that
\begin{equation}\label{eq.decomposition.of.V}
V=L_{0}\oplus L_0'=\langle e_{1},\dots,e_{n}\rangle\oplus \langle e^{1},\dots,e^{n}\rangle.
\end{equation}
Then any  $n$-plane  $L\in Gp_{n}(V)$ transversal to $L_{0}'$  has unique   basis $\{w_i\}$  projecting onto the basis $\{e_{i}\}$ (with respect to $L_{0}'$). Elements of such a basis can be written as
\begin{equation}\label{w_i}
w_{i}=e_{i}+p_{ij}e^{j},
\end{equation}
with the matrix $P=\left\Vert p_{ij}\right\Vert $ being symmetric if and only if $L$ is Lagrangian. So, every element $L\in\mathcal{L}(V)$ transversal to $L_{0}'$ is uniquely determined by a symmetric $n\times n$ real matrix $P$:
\begin{equation*}
L=L_P=\langle e_{i}+p_{ij}e^{j} \rangle
\end{equation*}
This gives  a local chart on $\mathcal{L}(V)$ with values in the vector space of symmetric matrices (hence, $\dim\mathcal{L}(V)=\frac{1}{2}n(n+1)$). It is
easy to check that coordinate changes in the overlaps between two such charts are $C^{\infty}$. The matrix $P$ of coordinates on $L$
transforms like a quadratic form
$$
P\mapsto\widetilde{P}=B^TPB
$$
where $B$ is the matrix of the transformation from basis $\{\widetilde{e}_i\}$ to basis $\{e_i\}$:
$\widetilde{e}_i\mapsto e_i=B^j_i\widetilde{e}_j$.

\smallskip\noindent
With respect to  a symplectic basis, an element of the symplectic
group $Sp(V)\simeq Sp_n(\mathbb{R})$ is represented by matrix
$$
\left(
\begin{array} {ll}
A & B \\
C & D
\end{array}
\right) \in Sp_n(\mathbb{R})
$$
with the blocks satisfying the conditions:%
\begin{equation*}
\left\{
\begin{array}
[c]{l}%
A^{T}\,C=C^{T}\,A\\
B^{T}\,D=D^{T}\,B\\
A^{T}\,D-C^{T}\,B=Id
\end{array}
\right.
\end{equation*}
The group  $Sp_n(\mathbb{R})$ acts  transitively on
$\mathcal{L}(V)$ by fractional linear transformations:
$$
Sp_n(\mathbb{R})\ni
\left(
\begin{array} {ll}
A & B \\
C & D
\end{array}
\right) : P\mapsto \widetilde{P}=(AP+B)(CP+D)^{-1}.
$$

\smallskip\noindent
We denote by $\mathcal{T}(\mathcal{L}(V))$ the {\bf tautological
bundle} of $\mathcal{L}(V)$, i.e. the vector bundle on
$\mathcal{L}(V)$ whose fibre at a point $L\in\mathcal{L}(V)$ is
the vector space $L$.

\smallskip\noindent
We have the {\bf Pl\"{u}cker embedding} of the Lagrangian
Grassmannian $\mathcal{L}(V)$ into the projective space
$\mathbb{P}\Lambda^n(V)$ given by
\begin{equation*}
\iota: L=\langle e_1,e_2,\dots , e_n \rangle \mapsto [\textrm{vol}_L]
\end{equation*}
where $\textrm{vol}_L=e_1 \wedge e_2 \wedge \cdots \wedge e_n$ is the volume element associated with the basis $\{e_i\}$ of $L$.

\smallskip\noindent
A straight line of the projective space $\mathbb{P}\Lambda^n(V)$
which is included in $\iota(\mathcal{L}(V))$ is called a {\bf line} of $\mathcal{L}(V)$. We will denote by $\ell(L,\dot{L})$ the line of $\mathbb{P}\Lambda^n(V)$ starting from $L$ in  direction
$\dot{L}\in T_L\mathcal{L}(V)$.

\smallskip\noindent
From now on, where needed, we shall identify $\mathcal{L}(V)$ with $\iota(\mathcal{L}(V))$.

\subsection{Metrics associated with tangent  and cotangent vectors of $\mathcal{L}(V)$}

Below we prove that the bundle $T\mathcal{L}(V)$ is canonically isomorphic to the symmetric square $S^2\big(\mathcal{T}^*(\mathcal{L}(V))\big)$ of the dual bundle $\mathcal{T}^*(\mathcal{L}(V))$ of the tautological bundle of $\mathcal{L}(V)$.\\
Namely, let $\dot{L}_{0} \in T_{L_{0}}\mathcal{L}(V)$ and
$\phi_{t}$  an  1-parameter  subgroup of $Sp(V)$ such that $\dot
L_0 =  \left.\frac{d\phi_t(L_0)}{dt}\right\vert_{t=0} $. The symmetric bilinear
form $g^{\dot L_0}$ on $L_0$ is defined by
\begin{equation}\label{g_Rdot}%
g^{\dot{L}_{0}}(v,w)\overset{\text{def}}{=}\omega\left(  \left.  \frac
{d\phi_{t}(v)}{dt}\right\vert _{t=0},w\right)  ,~~v,w\in L_{0}.
\end{equation}
It does not depend on 1-parametric  group $\phi_t$ whose orbit has
tangent vector  $\dot L_0$. Indeed,  any other  such 1-parameter
group can be written as $\phi'_t = \phi_t \circ h_t + o(t)$ where
$h_t $ belongs  to the  stabilizer $H = Sp(V)_{L_0}$ of  the point
$L_0$. Then
$$
\left.\frac{d\phi'_{t}(v)}{dt}\right\vert _{t=0}=
\left.\frac{d\phi_{t}(v)}{dt}\right\vert _{t=0} + \left.\frac{dh_{t}(v)}{dt}\right\vert
_{t=0}
$$
and
$$
\omega\left(  \left.  \frac {d\phi'_{t}(v)}{dt}\right\vert
_{t=0},w\right) = \omega\left(  \left.  \frac
{d\phi_{t}(v)}{dt}\right\vert _{t=0},w\right)
$$
since $\omega|_{L_0} =0$. Then we get the following theorem.
\begin{theorem}\label{th.vector=quadratic}
The map defined by \eqref{g_Rdot}
\begin{equation}\label{vector=quadratic}
g:T_{L}\mathcal{L}(V)  \longrightarrow  S^{2}(L^{\ast})\,,\quad
\dot{L}  \longmapsto  g^{\dot{L}}
\end{equation}
is a canonical isomorphism of the  tangent bundle $T\mathcal{L}(V)$ with the symmetric
square $S^2\big(\mathcal{T}^*(\mathcal{L}(V))\big)$ of the dual
tautological bundle.
\end{theorem}
In particular, a vector field $X$ on $\mathcal{L}(V)$ defines  a section $g^X$ of $S^2\big(\mathcal{T}^*(\mathcal{L}(V))\big)$ which we will call a metric on $\mathcal{T}(\mathcal{L}(V))$ (note that it can be \emph{degenerate}).\\
In terms of coordinates $p_{ij}$, the metric $g^{\dot L}$ on $L=\langle e_i + p_{ij}e^j \rangle$ associated with $\dot L \sim \dot P = ||\dot p_{ij}||$ is given by
$$
g^{\dot L} =  \dot p_{ij}e^i \otimes e^j.
$$
By duality, we get
\begin{corollary}
There is a canonical isomorphism
\begin{equation}\label{vector=quadratic.2}
g:T_{L}^\ast\mathcal{L}(V)  \longrightarrow  S^{2}(L)\,,\quad
\rho  \longmapsto  g_{\rho}
\end{equation}
of the cotangent bundle $T^*\mathcal{L}(V)$ with the symmetric square $S^2\big(\mathcal{T}(\mathcal{L}(V))\big)$ of the
tautological bundle.
\end{corollary}
There is no ambiguity in denoting by $g$ both the maps \eqref{vector=quadratic} and \eqref{vector=quadratic.2}: in fact vectors appear as superscripts whereas covectors as subscripts.\\
A 1-form $\rho $ on $\mathcal{L}(V)$ defines a section $g_{\rho}$ of $S^2\big(\mathcal{T}(\mathcal{L}(V))\big)$ which we call a metric on $\mathcal{T}^*(\mathcal{L}(V))$ (note that it can be \emph{degenerate}).

\smallskip\noindent
In terms of coordinates $p_{ij}$, the metric $g_{\rho}$ on
$L^*$ associated with $1$-form $\rho=\rho^{ij}dp_{ij}$, with $\|\rho^{ij}\|$ being the symmetric matrix of
coordinates of $\rho$ with respect to basis $\{(dp_{ij})_L\}$ of $T_{L}^{\ast}\mathcal{L}(V)$, is
\begin{equation}\label{eq.g.rho.local}
g_{\rho} = \rho^{ij}w_i \otimes w_j
\end{equation}
where $L = \langle w_i = e_i + p_{ij}e^j \rangle$. In particular, a function $F\in C^\infty(\mathcal{L}(V))$,
defines a metric on $L^*$:
\begin{equation}\label{eq.metric.dF.local}
g_{(dF)_L}=\sum_{i\leq j}\frac{\partial F}{\partial p_{ij}}w_i\vee w_j
\end{equation}
where we recall that $w_i\vee w_j=\frac{1}{2}(w_i\otimes w_j + w_j\otimes w_i)$.

\begin{remark}\label{rem.conformal.change}
Under conformal change $\omega \to \lambda \omega$ of  the
symplectic form, the above metrics change as
$$
g^{\dot L} \mapsto \lambda
g^{\dot L},\,\,   g_{\rho} \mapsto  \lambda^{-1} g_{\rho}.
$$
\end{remark}

\subsection{Lagrangian Grassmannian as  a homogeneous space}
The group $Sp(V)$  acts transitively on $\mathcal{L}(V)$  and the
stabilizer   $H$ of a point $L_0\in\mathcal{L}(V)$ is
$H=GL(L_0)\ltimes S^2(L_0)$. Hence we can identify
$\mathcal{L}(V)$ with the coset space
$$
\mathcal{L}(V)= Sp(V)/ \big(GL(L_0)\ltimes S^2(L_0)\big).
$$

\smallskip\noindent
Lagrangian Grassmanniann $\mathcal{L}(V)$ is a compact manifold and the maximal compact subgroup $U(n)$ of the group $Sp(V)=Sp(n,\mathbb{R})$ acts on it transitively with stabilizer $O(n)$. So we can identify
${\mathcal L}(V)$ with the symmetric space $U(n)/O(n)$, (whose  central symmetry at $o=e^{O(n)}$ is defined by complex conjugation). Note that the square of the determinant
$$
{\det}^2: U(n)/O(n) \to S^1
$$
defines a fibration over the circle $S^1$ with fibre $SU(n)/SO(n)$. The pull back $(\det^2)^*(d \varphi)$ of the  fundamental class $[d\varphi]$ of the circle is called the Maslov index of $\mathcal{L}(V)$.

\smallskip\noindent
The tautological bundle$ \mathcal{T}\mathcal{L}(V)$ is a
homogeneous vector bundle associated with the principal vector
bundle
$$
Sp(V)\to Sp(V)/H=\mathcal{L}(V)
$$
and the tautological representation
$$
H=GL(L_0)\ltimes S^2(L_0)\to GL(L_0)
$$
with kernel  $S^2(L_0)$.

\smallskip\noindent
Decomposition \eqref{eq.decomposition.of.V} induces  a gradation
of the Lie algebra $\mathfrak{sp}(V)$ of $Sp(V)$ (which is identified
with the symmetric square $S^2(V)$) given by
$$
\mathfrak{sp}(V)=\mathfrak{g}^{-1}+\mathfrak{g}^0+\mathfrak{g}^1 =
S^2(L'_0) + L'_0\vee L_0 + S^2(L_0).
$$
We  identify $\mathfrak{m}= \mathfrak{g}_{-1}=S^2(L'_0)$ with the
tangent space $T_{L_0}\mathcal{L}(V)$ and $\mathfrak{h}=L'_0\vee
L_0 + S^2(L_0)$ with the Lie algebra of the stabilizer $H$. The
commutative ideal $S^2(L_0)$ is the kernel of the isotropy
representation of $\mathfrak{h}$ on $\mathfrak{m}$ and the
stability  subalgebra $\mathfrak{h}= L_0'\vee L_0\simeq
\mathfrak{gl}(L_0')$ acts on $\mathfrak{m}$ in the natural way.
Hence we  get  an identification of  the tangent space
$T_{L_0}\mathcal{L}(V)$ with space of symmetric bilinear forms on
$L_0$ :
$$
T_{L_0}\mathcal{L}(V)  \simeq S^2(L_0^*).
$$
According to Theorem
\ref{th.vector=quadratic}, this identification does not depend on
the choice of $L_0'$.
\smallskip\noindent
Note that in terms of basis $\{e_i\}$ of $L_0$ and the dual basis
$\{e^i\}$ of $L_0'\simeq L_0^*$, the matrix of elements of
$\mathfrak{sp}(V)$ has the form
$$
\left(
\begin{array} {cc}
A & B \\
C & -A^T
\end{array}
\right)
$$
where $ A \in \mathfrak{gl}(L_0),\, B \in S^2(L_0), \, C \in
S^2(L_0')$.

\subsection{Rank of tangent vectors of $\mathcal{L}(V)$ and its geometrical meaning}

By using Theorem \ref{th.vector=quadratic}, we define the {\bf
rank} of a tangent vector $\dot{L}\in T\mathcal{L}(V)$ as the rank of the
corresponding bilinear symmetric forms $g^{\dot L}$. In view of Remark \ref{rem.conformal.change}, this definition is invariant under a conformal change of the symplectic form. Of course, proportional tangent vectors have the same rank. We denote by
$$
T^k\mathcal{L}(V)=\{\dot{L}\in T\mathcal{L}(V) \,\,|\,\, \text{rank}(\dot{L})=k\}
$$
the set of vectors of rank $k$ and define the canonical map $\Rad: T\mathcal{L}(V)\to Gr_{n-k}(V)$ which associates with any tangent vector $\dot{L}\in T\mathcal{L}(V)$ the radical of $g^{\dot{L}}$:
\begin{equation}\label{eq.Rad}
\Rad(\dot{L}):=\Rad\,(g^{\dot{L}}).
\end{equation}
In the next section we shall construct a sort of inverse of  map $\Rad$ (see Remark \ref{rem.rank1.and.hyperplane}). Now we give a geometrical interpretation of $\Rad(\dot{L})$. The space $\Rad(\dot{L})$ is the intersection of the plane
$L$ and the infinitesimally close Lagrangian plane $L+\dot{L}dt$, more precisely,
$$
\Rad(\dot{L})=\lim_{t\to 0} \, L\cap L(t)\,,\,\,\, L(0)=L, \,\, \dot{L}(0)=\dot{L}.
$$
Indeed if
$
L=\{ x=x^ie_i \}
$
and
$
L(t)= \{ x^i(e_i + p_{ij}(t)e^j) \}
$
then
$$
L\cap L(t) =\{ x=x^ie_i\,\,|\,\, p_{ij}(t)x^i=0\}=\Rad(P(t))
$$
and $\Rad(\dot{L})=\lim_{t\to 0}\, L\cap L(t) = \Rad(\dot{P}(0))$.

\smallskip\noindent
We call the set $T^1\mathcal{L}(V)$ of vectors of rank $1$  the  \textbf{characteristic cone} or \textbf{Segre variety} (see \cite{AG}). If $\dot{L}\in T^1\mathcal{L}(V)$, then, up to a sign,
\begin{equation}\label{eq.rank.1.vectors}
\dot{L}\simeq g^{\dot{L}}= \eta\otimes\eta\,,\,\,\,\text{for some }\eta\in L^*
\end{equation}
and the canonical map $\Rad$ takes values in $Gr_{n-1}(L)\simeq\mathbb{P}L^*$. From now on, unless otherwise specified, we identify $\dot{L}$ with $g^{\dot{L}}$.

\smallskip\noindent
In terms of coordinates, if $ L=\langle w_i=e_i+p_{ij}e^j \rangle
$ and $\dot{L}\in T^1\mathcal{L}(V)$ has coordinates
$\dot{p}_{ij}$, then by \eqref{eq.rank.1.vectors}
$\dot{p}_{ij}=\eta_i\eta_j$ and
\begin{equation*}
\Rad(\dot{L})= [\eta_ie^i]\in \mathbb{P}L^*.
\end{equation*}

\smallskip\noindent
We recall the straight line $\ell(L,\dot{L})$ in
$\mathbb{P}\Lambda^n(V)$ starting from $L$ in direction
$\dot{L}\in T_L\mathcal{L}(V)$.
\begin{proposition}\label{prop.line}
The straight line  $\ell(L, \dot L)$ of $\mathbb{P}\Lambda^n(V)$ is a line of $\mathcal{L}(V)$ (i.e. it is included in $\mathcal{L}(V)$) if and only if
$\text{rank}(\dot L) =1$, i.e. $\dot L \in T^1_L\mathcal{L}(V)$.
\end{proposition}
To prove the proposition we need the following lemma.
\begin{lemma}
Let $a,a'\in\Lambda^k(W)$ be two $k$-vectors such that $ta+sa'$ is decomposable for any $t,s\in\mathbb{R}$. Then there exists a decomposable $(k-1)$-vector
$b\in\Lambda^{k-1}(W)$ and vectors $v$, $v'$ such that
$a=v\wedge b$ and $a'=v'\wedge b$.
\end{lemma}
\begin{proof}
A $k$-vector $c$ is decomposable iff it satisfies  the Pl\"{u}ker
relation  $(\gamma\,\lrcorner\, c)\wedge c=0$ for any $\gamma\in
\Lambda^{k-1}(W^*)$ (see, for example \cite{GH}). By hypothesis
these relations hold for $c=a$, $c=a'$ and $c=a+a'$. Then we
derive that
$$
0=(\gamma\,\lrcorner\,a)\wedge a' + (\gamma\,\lrcorner\,a')\wedge a\,, \,\,\,\forall\,\,\gamma\in \Lambda^{k-1}(W^*).
$$
We choose $\gamma$ such that $v':=\gamma\,\lrcorner\,a\neq 0$ and
$v:=-\gamma\,\lrcorner\,a'\neq 0$. Then $ v'\wedge a = v\wedge a',
$ so that $ a = v\wedge b\,,\,\,\,a' = v'\wedge b $ for some
$b\in\Lambda^{k-1}(W)$.
\end{proof}

\smallskip
\begin{proof}[Proof of Proposition \ref{prop.line}]
Assume that $\dot L \in T^1_L\mathcal{L}(V)$. We can choose
local coordinates $P = ||p_{ij}||$ such that $ P(L)=0$ and
 $P(\dot L) =\text{diag}(1,0,\dots,0)$. Then the straight line
$$
\ell(L, \dot L) = [(e_1 + t e^1)\wedge e_2 \cdots \wedge e_n]=
[e_1 \wedge \cdots \wedge e_n + te^1\wedge e_2\wedge\dots\wedge e_n]
$$
is included in $\mathcal{L}(V)$.

\smallskip\noindent
The converse claim  follows from the above lemma.
\end{proof}

\section{Submanifolds of the Lagrangian Grassmannian $\mathcal{L}(V)$}\label{sec.sub.of.Lagr.Grass}

\subsection{Characteristic cone and
characteristic subspaces of a hypersurface $\mathrm{E}$ of $\mathcal{L}(V)$ and its conformal metric $g_{\mathrm{E}}$}\label{sec.char.cone.Lagr.Grass}

Let
$$
\mathrm{E}=\{F=0\}
$$
be a hypersurface of $\mathcal{L}(V)$ which is the zero level set of a non singular function $F\in C^{\infty}(\mathcal{L}(V))$. We denote  by
$$
g_{\mathrm{E}}:=[g_{dF}|_{\mathrm{E}}]\,,
$$
the conformal class of the restriction to $\mathrm{E}$ of the contravariant metric $g_{dF}$. It is easy to se that $g_{\textrm{E}}$ depends only on the
hypersurface $\textrm{E}$ and is called the {\bf conformal metric} associated with $\textrm{E}$. Its local expression is given by \eqref{eq.metric.dF.local}.

\begin{definition}\label{def.cono.car.point}
The set
$$
\text{Ch}_{L}(\mathrm{E})=
T_{L}\mathrm{E}\cap T^{1}_{L}\mathcal{L}(V)
$$
of rank $1$ tangent vectors to $\mathrm{E}$
is called the \textbf{characteristic cone} at $L$ of the hypersurface $\mathrm{E}$. Elements of $\text{Ch}_{L}(\mathrm{E})$ are called \textbf{characteristic vectors} for $\mathrm{E}$ at $L$. The $1$-dimensional vector space generated by a characteristic vector is called a \textbf{characteristic direction}. A characteristic vector $\dot{L}$ for $\mathrm{E}$ at $L$ is called \textbf{strongly characteristic} if the associated line $\ell(L,\dot{L})$ is contained in $\mathrm{E}$.
\end{definition}

\begin{proposition}\label{prop.char.isotropic}
Characteristic vectors $\dot L \in \text{Ch}_{L}(\mathrm{E})$  are, up to sign, the tensor square $\dot L =  \eta \otimes \eta$ of $g_{\mathrm{E}}$-isotropic  covectors $\eta \in  L^*$.
\end{proposition}
\begin{proof}
A tangent vector $\dot{L} \in T_L\mathcal{L}(V)$  with coordinates
$\dot{P}=||\dot{p}_{ij}||$  has rank $1$ iff $\dot{p}_{ij}= \pm\eta_i\eta_j$ (see \eqref{eq.rank.1.vectors}). It is characteristic for $\mathrm{E}$ at $L$ if and only if
\begin{equation}\label{eq.characteristics}
\sum_{i\leq j}\frac{\partial F}{\partial{p_{ij}}}\,\dot{p}_{ij}=
\sum_{i\leq j}\frac{\partial
F}{\partial{p_{ij}}}\,\eta_i\eta_j=g_{\mathrm{E}}(\eta,\eta)=0\,,
\end{equation}
i.e. iff the covector $\eta=\Rad(\dot{L})$  is
$g_{\mathrm{E}}$-isotropic.
\end{proof}

\medskip\noindent
We define the \textbf{prolongation}
$U^{(1)}\subset\mathcal{L}(V)$ of a subspace $U \subset V$ by
:
\begin{equation}\label{eq.prolongation.of.U}
U^{(1)}:=\left\{
\begin{array}{c}
L\in \mathcal{L}(V)\,\,|\,\, L\supseteq U,\,\,\text{if}\,\,\dim(U)\leq n\\
\\
L\in \mathcal{L}(V)\,\,|\,\, L\subseteq U,\,\,\text{if}\,\,\dim(U)\geq n\\
\end{array}
\right.
\end{equation}
Since $L=L^\perp$, one  can easily check that
\begin{itemize}
\item $U\subset W \Longrightarrow U^{(1)}\supset W^{(1)}$;
\item $U^{(1)}=\left({U^\perp}\right)^{(1)}$.
\end{itemize}
The following simple proposition describes  the prolongation
$U^{(1)}$ of an isotropic subspace $U$ of $V$.
\begin{proposition}\label{prop.prol.U}
Let $U$ be an isotropic $k$-dimensional subspace of $V$. Let $U'$
be also an isotropic $k$-dimensional subspace of $V$ such that
$\omega$ is not degenerate on $U\oplus U'$. Then $W:=(U\oplus
U')^\perp$ is a symplectic subspace  and
$$
U^{(1)}\simeq U\oplus \mathcal{L}(W):=\{U\oplus L'\,|\, L'\in \mathcal{L}(W)\}.
$$
In particular
\begin{equation}\label{eq.dim.prol.U}
\dim\,U^{(1)}=\dim\mathcal{L}(W)=\frac{(n-k)(n-k+1)}{2}.
\end{equation}
\end{proposition}

\begin{definition}\label{def.char.point}
An isotropic subspace $U$ is called \textbf{characteristic} for a covector $\rho\in T^*_L\mathcal{L}(V)$ if $U\subset L$ and
$\rho|_{T_LU^{(1)}}=0$. It is called characteristic for a
hypersurface $\mathrm{E}=\{F=0\}$ of $\mathcal{L}(V)$ at a point
$L\in\mathrm{E}$ if it is characteristic for $(dF)_L$. It is
called \textbf{strongly characteristic} if
$U^{(1)}\subset\mathrm{E}$. A covector $\eta\in L^*$ is called characteristic for $\rho$ if $\Ker(\eta)$ is characteristic for $\rho$.
\end{definition}
\begin{remark}\label{rem.char}
Previous definition is also valid for submanifolds of $\mathcal{L}(V)$ of any dimension. We restrict our attention to hypersurfaces of $\mathcal{L}(V)$ as our target is to treat characteristics of scalar second order PDEs with one unknown function (see Section \ref{sec.char}).
\end{remark}
The following remark clarifies the relationship between characteristic directions and characteristic subspaces.
\begin{remark}\label{rem.rank1.and.hyperplane}
Prolongation \eqref{eq.prolongation.of.U} is a sort of inverse of map \eqref{eq.Rad}. Namely, any $\dot L = \pm \eta \otimes \eta
\in T^1_L\mathcal{L}(V)$ defines the hyperplane $H=\Rad(\dot{L})=\Ker(\eta)$ of $L$ which has the property that $T_LH^{(1)}=\langle \dot L \rangle$, and
viceversa (we note that $H^{(1)}$ is $1$-dimensional in view of \eqref{eq.dim.prol.U}). Thus we have the following correspondence:
\begin{equation}\label{eq.correspondence.gianni}
\text{hyperplanes of $L$ (which correspond  to  elements of
$\mathbb{P}L^*$)$\,\,\,\Longleftrightarrow\,\,\,$ directions of
$T_L\mathcal{L}(V)$ of rank $1$}
\end{equation}
It follows that if $\Ker(\eta)=H \subset L$ is a hyperplane of a Lagrangian plane $L$ then $H^{(1)}=\ell(L,\dot{L}=\eta\otimes\eta)=\{L_t\}$
is a straight line of $\mathcal{L}(V)$ in view of Proposition
\ref{prop.line}. Restricting \eqref{eq.correspondence.gianni} to a hypersurface $\mathrm{E}$ of $\mathcal{L}(V)$ we have the following correspondence:
\begin{equation}\label{eq.correspondence.gianni.2}
\text{$(n-1)$-dimensional characteristic subspaces for $\mathrm{E}$ at $L$  $\,\,\,\Longleftrightarrow\,\,\,$ characteristic directions for $E$ at $L$}
\end{equation}
\end{remark}
We have already seen, in Proposition \ref{prop.char.isotropic}, that a vector $\dot{L}= \pm \eta \otimes \eta
\in T_L\mathcal{L}(V)$ is characteristic for $\mathrm{E}$ at $L$ if $ \eta \in L^*$ is  $g_{\mathrm{E}}$-isotropic. Next theorem generalizes this property.
\begin{theorem}\label{th.fab}
Let $U\subset L\in\mathcal{L}(V)$ and $\rho\in
T_{L}^{\ast}\mathcal{L}(V)$. Then $U$ is characteristic for
$\rho$ if and only if its annihilator $ U^0 \subset L^{\ast}$ is
$g_{\rho}$-isotropic.
\end{theorem}
\begin{proof}
Let $\{e_{1},\dots,e_{n}\}$ be a basis of $L$ such that $\{e_{a}\}_{a = 1, \cdots, k}$ is a basis of $U$. Let also $\{e_{1},\dots,e_{n},e^{1},\dots,e^{n}\}$ be its extension to a symplectic basis of $V$.
Then we can consider $\{e^i\}_{i = k+1, \cdots , n}$ as a
basis of $U^0$. So $U^0$ is $g_\rho$-isotropic if
$$
g_{\rho}(e^i, e^j) = \rho^{ij} =0, \,\, i,j \in \{k+1,\cdots,n \},
$$
with $g_{\rho}$ as in \eqref{eq.g.rho.local}. By Proposition \ref{prop.prol.U},
$$
U^{(1)} = \{ L = \langle e_a, e_i + p_{ij}e^j \rangle\,\,\,\big|\,\,\, 1\leq a\leq k\,\,,\,\,||p_{ij}|| \in S^2\mathbb{R}^{n-k}  \}.
$$
Then its tangent space is given by
$$
T_{L}U^{(1)}=\langle e^i \vee e^j
 ,\,\, i,j=k+1,\dots,n\rangle.
$$
Hence, $U$ is characteristic for $\rho$ if and only if%
\[
\rho( e^i \vee e^j)= \rho^{ij}=0,   ~i,j = k+1 ,\dots,n
\]
which means  that $U^0$ is $g_\rho$-isotropic.
\end{proof}

\begin{corollary}
Let $F = F(p_{ij})$ be a function on $\mathcal{L}(V)$. Then a subspace $U \subset L$, in view of \eqref{eq.metric.dF.local}, is characteristic for $(dF)_L$ (i.e. for the hypersurface $\mathrm{E}=\{F=0\}$ at $L$) iff
$$
g_{(dF)_L}(\alpha,\beta)=\frac{1}{2}\sum_{i\leq j}\frac{\partial
F}{\partial{p_{ij}}}(\alpha_i\beta_j+\alpha_j\beta_i)=0\,, \,\,\, \forall \,
\alpha,\beta \in U^0.
$$
\end{corollary}
In view of previous theorem we have the following correspondence:
$$
\text{$\eta$ is characteristic for $\rho\,\, \Longleftrightarrow\,\,\eta$ is $g_{\rho}$-isotropic $\,\, \Longleftrightarrow\,\,\rho(\eta\otimes\eta)=0$}.
$$
In the case in which $\rho=dF$, the last property means that the vector $\eta\otimes\eta$ is characteristic for $\{F=0\}$ at the point $L$ (see also Remark \ref{rem.rank1.and.hyperplane}).

\begin{theorem}\label{th.rotating.debole}
Let $\rho\in T^*_L\mathcal{L}(V)$. Then $g_\rho$ is decomposable
iff $(n-1)$-dimensional characteristic subspaces for $\rho$ (at
$L$) form two $(n-2)$-parametric families $\mathcal{H}$ and
$\mathcal{H}'$ such that
$$
\dim\bigcap_{U\in \mathcal{H}}U = \dim\bigcap_{U\in \mathcal{H}'}U =1.
$$
\end{theorem}
\begin{proof}
Assume that  $g_\rho$ is decomposable, i.e.  $g_\rho=v\vee w$ for some $v,w\in L$. By Theorem \ref{th.fab}, a
hyperplane $U=\Ker(\alpha)$ of $L$ is characteristic iff
$g_\rho(\alpha,\alpha)=\alpha(v)\alpha(w)=0$. This means that
$v\in U$ or $w\in U$. So we get two families
$$
\mathcal{H}=\{U\subset L\,\,|\,\,\ v\in U\}\,, \quad \mathcal{H}'=\{U\subset L\,\,|\,\, w\in U\},
$$
 of characteristic hyperplanes such that
$$
\bigcap_{U\in\mathcal{H}}U=\langle v \rangle\,, \quad \bigcap_{U\in\mathcal{H}'}U=\langle w \rangle.
$$
.\\
Assume now that  $\mathcal{H}$ is  a $(n-2)$-parametric family of
characteristic hyperplanes of $L$ which contain a common line
$\langle v \rangle$. By dimensional reason, the set
$$
\bigcup_{U\in\mathcal{H}} U^0=\{\alpha\in L^*\,\,|\,\,
\alpha|_U=0\,\,\,\text{for some $U\in\mathcal{H}$}\}
$$
contains a conic convex open subset $\mathcal{O}$ of the
annihilator $v^0 \subset L^*$. So $\alpha,\alpha'\in\mathcal{O}$
implies that $\alpha+\alpha'\in\mathcal{O}$. Theorem \ref{th.fab} shows that
$$
g_\rho(\alpha,\alpha)=g_\rho(\alpha',\alpha')=g_\rho(\alpha+\alpha',\alpha+\alpha')=0
$$
which  implies $g_\rho(\alpha,\alpha')=0$. Hence  any linear
combination of covectors in $\mathcal{O}$ is $g_\rho$-isotropic.
The set of such linear combinations coincides with the annihilator
$v^0$. The $g_\rho$-isotropy of all vectors in $v^0$ implies
that  $v^o$  is $g_\rho$-isotropic. Then
$$
g_\rho=v\vee w
$$
for some vector $w\in L$.
\end{proof}
\begin{remark}
The second part of the above proof shows that the existence of only one of the families of Theorem \ref{th.rotating.debole} implies the existence of the other one. Also, as by-product, we derive that each of such family consists of \emph{all} hyperplanes of $L$ containing some line.
\end{remark}

\subsection{Hypersurfaces $\mathrm{E}_\Omega$ of $\mathcal{L}(V)$ associated with $n$-forms $\Omega$ on $V$
and their characteristics}\label{sec.hyper.Lagr.Grass.assoc.n.forms}

Any $n$-form $\Omega\in \Lambda^n(V^*)$ defines the hypersurface
\begin{equation}\label{eq.E.Omega}
\mathrm{E}_\Omega=\{L\in\mathcal{L}(V)\,\,\,|\,\,\, \Omega|_L=0\}.
\end{equation}
That $\mathrm{E}_{\Omega}$ has codimension $1$ follows from the fact that, if $P=\|p_{ij}\|$ is the local chart on $\mathcal{L}(V)$ defined as in Section \ref{sec.Lagr.Grass}, then
$$
\Omega|_{L_P}=F(P)\,e_1^*\wedge\cdots\wedge e_n^* \simeq F(P)\,e^1\wedge\cdots\wedge e^n,
$$
for some function $F\in C^\infty(\mathcal{L}(V))$, $\{e^*_i\}$ being the dual basis of $\{w_i\}$ defined by
\eqref{w_i}; so, the condition in \eqref{eq.E.Omega} reduces to the vanishing of $F$.
\begin{remark}\label{rem.horizontalization.point}
The correspondence $L\in\mathcal{L}(V)\mapsto
\Omega|_L\in\Lambda^n(L^*)$ defines an $n$-form on the
tautological bundle $\mathcal{T}(\mathcal{L}(V))$ of $\mathcal{L}(V)$.
\end{remark}
Two $n$-forms $\Omega$, $\widetilde{\Omega}$ define the same
hypersurface ($\mathrm{E}_\Omega=\mathrm{E}_{\widetilde \Omega})$
if, up to a non vanishing factor, they are related by
\begin{equation*}
\widetilde{\Omega}=\Omega+\sigma\wedge\omega=\colon\Omega^\sigma
\end{equation*}
for some $\sigma\in \Lambda^{n-2}(V^*)$.

\smallskip\noindent
Note that hypersurfaces of the form $\mathrm{E}_{\Omega}$ can be obtained as intersections of $\mathcal{L}(V)$ (or, rather, its Pl\"{u}cker image) with hyperplanes of $\mathbb{P}\Lambda^n(V)$. In fact, such hyperplanes biunivocally correspond to hyperplanes of $\Lambda^n(V)$, which in their turn can be identified with lines in $\Lambda^n(V)^*$:
$$
(\mathbb{P}\Lambda^n(V))^*\simeq\mathbb{P}(\Lambda^n(V)^*);
$$
on the other hand, one can associate with any $\Omega\in\Lambda^n(V^*)$ the covector $\widetilde{\Omega}\in\Lambda^n(V)^*$ given by
$$
\widetilde{\Omega}(v_1\wedge\cdots\wedge v_n):=\Omega(v_1,\dots,v_n)\,,\,\,\, v_1,\dots,v_n\in V,
$$
so that $\Lambda^n(V^*)$ is canonically isomorphic to $\Lambda^n(V)^*$. Therefore,
$$
\mathrm{E}_{\Omega}=\mathcal{L}(V)\cap \{L\in\mathcal{L}(V)\subset\mathbb{P}\Lambda^n(V)\,\,|\,\,\widetilde{\Omega}(L)=0\}.
$$

\begin{theorem}\label{th.strongly.char}
Let $L\in\mathrm{E}_{\Omega}$. If a hyperplane $H$ of $L$ is
characteristic for $\mathrm{E}_{\Omega}$ at $L$ then it is
strongly characteristic.
\end{theorem}
\begin{proof}
Let us choose a symplectic basis $\{ e_i,e^i \}$ of $V$ such that
$H=\langle e_1,\dots,e_{n-1} \rangle$. Then $H^{(1)}=\{L_t=\langle
e_1,\dots,e_{n-1},e_n+te^n \rangle\}$. All Lagrangian planes in a neighborhood of $L$ are described by
$$\widetilde L=\langle e_i+p_{ij}e^j\rangle.$$
So we can define
$$
\mathrm{vol}_{\widetilde L}:=(e_1+p_{1j}e^j)\wedge\cdots\wedge(e_n+p_{nj}e^j).
$$
Also, for short, $\mathrm{vol}_t:=\mathrm{vol}_{L_t}$.

\smallskip\noindent
If $L'=\langle e_1,\dots,e_{n-1},e^n \rangle$, we add the notation
$\mathrm{vol}_{L'}:=e_1\wedge \cdots e_{n-1}\wedge e^n$ in such a way that
$\textrm{vol}_t=\textrm{vol}_L+t\textrm{vol}_{L'}$. In this way
the tangent vector to $H^{(1)}$ at $L$ is defined by the
derivative along $\textrm{vol}_{L'}$. Also, we define
$F(\widetilde{L})=\textrm{vol}_{\widetilde{L}}\,\lrcorner\,\Omega$ so that $\mathrm{E}_{\Omega}$ is locally described by $\{F=0\}$. The derivative of $F$ at
$L$ along $\textrm{vol}_{L'}$ is
\begin{multline*}
\lim_{t\to 0} \frac{F(L_t)-F(L)}{t} = \lim_{t\to 0}
\frac{\textrm{vol}_t\,\lrcorner\,\Omega-\textrm{vol}_L\,\lrcorner\,\Omega}{t}
= \lim_{t\to 0}
\frac{(\textrm{vol}_L+t\,\textrm{vol}_{L'})\,\lrcorner\,\Omega-\textrm{vol}_L\,\lrcorner\,\Omega}{t}=
\\
=
\textrm{vol}_{L'}\,\lrcorner\,\Omega=\Omega(e_1,\dots,e_{n-1},e^n)
\end{multline*}
which vanishes if and only if $L'$ belongs to
$\mathrm{E}_{\Omega}$. In this case we derive that $H^{(1)}$ is
included in $\mathrm{E}_{\Omega}$.
\end{proof}

\smallskip\noindent
Below we describe $(n-1)$-dimensional characteristic subspaces for
the hypersurface $\mathrm{E}_\Omega$ . We need the following
definition.
\begin{definition}
Let $\Omega\in\Lambda^n(V^*)$ be an $n$-form on a vector space
$V$. A $k$-dimensional subspace $U =\langle e_1, \cdots , e_k\rangle \subset
V$ is called \textbf{$\Omega$-isotropic} if $ (e_1 \wedge \cdots \wedge e_k
)\lrcorner\,\Omega=0$.
\end{definition}
%
Note that an $n$-dimensional subspace $U $ is $\Omega$-isotropic
if $\Omega|_U=0$. Next theorem describes $(n-1)$-dimensional
characteristic subspaces of $\mathrm{E}_{\Omega}$.
\begin{theorem}\label{th.car.hyp}
Let $L\in\mathrm{E}_{\Omega}$. A hyperplane $H$ of $L$ is
characteristic for $\mathrm{E}_{\Omega}$ at $L$ iff $H$ is $\Omega^\sigma$-isotropic for some $\sigma\in\Lambda^{n-2}(V^*)$.
\end{theorem}
\begin{proof}
We use  the same notations as in the proof of Theorem
\ref{th.strongly.char}.

\smallskip\noindent
Then
$$
H \text{ is characteristic }\, \Longleftrightarrow \,
H^{(1)}\subset\mathrm{E}_{\Omega}\, \Longleftrightarrow\,
\textrm{vol}_t\,\lrcorner\,\Omega=0 \, \Longleftrightarrow
\Omega_a(e_n)=\Omega_a(e^n)=0
$$
where $
\Omega_a=a\,\lrcorner\,\Omega\,,\,\,\,a=e_1\wedge\cdots\wedge
e_{n-1}. $ For any $\sigma\in\Lambda^{n-2}(V^*)$, we have that
$$
a\,\lrcorner\,\Omega^\sigma=\Omega_a + \sum_{j} (-1)^j
\sigma(e_1,\dots,e_{j-1},e_{j+1},\dots,e_{n-1})(e_j\,\lrcorner\,\omega).
$$
In particular, $(a\,\lrcorner\,\Omega^\sigma)|_{L'}=0$ and
$$
(a\,\lrcorner\,\Omega^\sigma)(e^i)=\Omega_a(e^i) + (-1)^i
\sigma(e_1,\dots,e_{i-1},e_{i+1},\dots,e_{n-1})
$$
which vanishes if
$$
\sigma(e_1,\dots,e_{i-1},e_{i+1},\dots,e_{n-1})=(-1)^{i+1}\Omega_a(e^i).
$$
Then, for such $\sigma$, $a\,\lrcorner\,\Omega^\sigma=0$, i.e. $H$
is isotropic for $\Omega^\sigma$.

\smallskip\noindent
The converse statement is trivial. In fact, if $H=\langle
e_1,\dots,e_{n-1}\rangle$ is $\Omega$-isotropic, then $\Omega_a=0$
which implies
$e_n\,\lrcorner\,\Omega_a=e^n\,\lrcorner\,\Omega_a=0$.
\end{proof}

\begin{remark}
If $H$ is an isotropic $(n-1)$-plane which contains at least one vector of $\Ker\,\Omega^\sigma$ then it is
$\Omega^\sigma$-isotropic and hence characteristic. Converse statement is not true: it may happen that a characteristic plane $H$ has trivial intersection with the kernels of all forms of type $\Omega^\sigma$, $\sigma\in\Lambda^{n-2}(V^*)$. For instance, for $n=3$, consider the following example:
$$
H=\langle e_1,e_2 \rangle,\quad \Omega=e_1^*\wedge e_3^*\wedge
e^{2*}+ e^{2*}\wedge e^{1*}\wedge e^{3*},
$$
where $\{e_i,e^i\}$ is a symplectic basis. However the following proposition says that this is true for
decomposable $n$-forms.
\end{remark}
\begin{proposition}\label{prop.appoggio100}
An  $(n-1)$-dimensional  subspace $H$  is $\Omega$-isotropic  for a  decomposable $n$-form  $\Omega$ if and only if
$$
H\cap\Ker\,\Omega\ne 0.
$$
\end{proposition}
\begin{proof}
Let $\Omega=\varrho_1\wedge\cdots\wedge\varrho_n$ and $H=\langle e_1,\dots, e_{n-1}\rangle$ such that
$\textrm{vol}_H\,\lrcorner\,\Omega=0$. It implies that rank of the
$||\varrho_i(e_j)||$ is $\leq n-2$. Hence there exists a linear
combination $e:=\lambda^je_j\in H$ such that $\varrho_i(e)=0$,
which entails $e\in\Ker\,\Omega$.
\end{proof}

\subsection{Hypersurfaces $\mathrm{E}_D$ of $\mathcal{L}(V)$ associated with an $n$-plane $D\subset V$ and their characteristics}\label{sec.variety.ED.point}

\subsubsection{Definition of $\mathrm{E}_D$ and reconstruction of $D$ from $\mathrm{E}_D$  }

We associate with an $n$-dimensional subspace $D\subset V$ the subset of $\mathcal{L}(V)$
$$
\mathrm{E}_D=\{ L\in\mathcal{L}(V)\,\,|\,\,L\cap D\neq 0 \}
$$
consisting of all Lagrangian planes which non trivially
intersect $D$. With respect to a symplectic basis $\{e_i,e^i\}$  the subspace $D$ can be written
as
\begin{equation}\label{eq.appoggio.D}
D=\langle w_i=e_i+b_{ij}e^j \rangle = \{x=x^ie_i+x^ib_{ij}e^j\}
\end{equation}
where $B=||b_{ij}||$ is an $n\times n$ matrix.
If we denote by $D^\perp$ the \textbf{orthogonal complement} of $D$ w.r.t. the symplectic form $\omega$, we have that
\begin{equation*}
D^\perp=\langle w'_i=e_i+b_{ji}e^j \rangle.
\end{equation*}
In particular, $D$ is a Lagrangian plane iff matrix $B$ is symmetric, as $D=D^\perp$. The proposition below shows that $\mathrm{E}_D$ is an algebraic hypersurface of $\mathcal{L}(V)$.
\begin{proposition}\label{prop.local.ED.point}
In terms of the coordinates $P=||p_{ij}||$ of $L= L_P \in
\mathcal{L}(V)$ associated with the basis $\{e_i,e^i\}$, $\mathrm{E}_D$ is described as follows:
$$
\mathrm{E}_D=\{L_P\,\,\,\big|\,\,\, \det(P-B)=0\}
$$
with $D$ given by \eqref{eq.appoggio.D}.
\end{proposition}
\begin{proof}
Since
$$
L=\langle e_i+p_{ij}e^j \rangle = \{x=x^ie_i+x^ib_{ij}e^j\}.
$$
we have
$L\cap D=\{x=x^ie_i+x^ib_{ij}e^j\,\,|\,\,(P-B)\cdot
x=0\} =\Ker\, (P-B). $
\end{proof}

\smallskip\noindent
Equations of type $\mathrm{E}_D$ are also defined by $n$-forms (and then are of the type introduced in Section \ref{sec.hyper.Lagr.Grass.assoc.n.forms}) as the following proposition shows.
\begin{proposition}\label{prop.n.forms}
Let $D=\{\varrho_1=\varrho_2=\dots=\varrho_n=0\}$ be an
$n$-dimensional subspace defined by $n$ linear forms, then
$\mathrm{E}_D = \mathrm{E}_{\Omega_D}$ where
\begin{equation*}
{\Omega_D}=\varrho_1\wedge\cdots\wedge\varrho_n.
\end{equation*}
\end{proposition}

\begin{theorem}\label{th.D.equal.D.orth}
Let $(V,\omega)$ be a $2n$-dimensional symplectic vector space. Let $D$ and $\widetilde{D}$ be $n$-dimensional planes of $V$. Then
$$
\mathrm{E}_{\widetilde{D}}=\mathrm{E}_{D}\,\,\Longleftrightarrow  \widetilde{D}=D \text{\,\,\,or\,\,\,} \widetilde{D}=D^\perp.
$$
\end{theorem}
\begin{proof}
The condition is necessary. Let $e\in\widetilde{D}$ so that $e^{(1)}\subset \mathrm{E}_{\widetilde{D}}=\mathrm{E}_{D}$ which implies that $L\cap D\neq 0$ for all $L\in e^{(1)}$. We shall prove that $e \in D$ or $e \in D^{\perp}$. Choose  a symplectic basis $\{e_i,e^i\}$ such that $e_1=e$ and
$$
D=\langle e_i+ b_{ij} e^j \rangle,
$$
for some $b_{ij}\in\mathbb{R}$.

\smallskip\noindent
Then the vector $e=e_1$ belongs to $D$ iff $b_{1j}=0$ for any $j$ and belongs to $D^\perp$ iff $b_{j1}=0$. We shall show that if all Lagrangian subspaces containing the vector $e$ intersects $D$ (non trivially), then either $b_{1j}=0$ or $b_{j1}=0$.

\smallskip\noindent
In order to do this, we shall choose appropriated Lagrangian
subspaces.

\smallskip\noindent
Let us consider the Lagrangian subspace
\begin{equation*}
L=\langle e_1,e^2,\dots,e^n \rangle.
\end{equation*}
By hypothesis $L$ intersects non trivially $D$. So
the determinant of the following matrix
\begin{equation*}
\left(
\begin{array}{ccccc|ccccc}
1 & 0 & 0 & \cdots  & 0 & 0 & 0 & 0 & \cdots  & 0 \\
0 & 0 & 0 & \cdots  & 0 & 0 & 1 & 0 & \cdots  & 0 \\
0 & 0 & 0 & \cdots  & 0 & 0 & 0 & 1 & \cdots  & 0 \\
\vdots  & \vdots  & \vdots  & \ddots  & \vdots  & \vdots  & \vdots  & \vdots
& \ddots  & \vdots  \\
0 & 0 & 0 & \cdots  & 0 & 0 & 0 & 0 & \cdots  & 1 \\
\hline
1 & 0 & 0 & \cdots  & 0 & b_{11} & b_{12} & b_{13} & \cdots  & b_{1n} \\
0 & 1 & 0 & \cdots  & 0 & b_{21} & b_{22} & b_{23} & \cdots  & b_{2n} \\
\vdots  & \vdots  & \ddots  & \ddots  & \vdots  & \vdots  & \vdots  & \vdots
& \ddots  & \vdots  \\
0 & 0 & 0 & \ddots  & 0 & b_{n-1\;1} & b_{n-1\;2} & b_{n-1\;3} & \cdots  &
b_{n-1\;n} \\
0 & 0 & 0 & \cdots  & 1 & b_{n1} & b_{n2} & b_{n3} & \cdots  & b_{nn}
\end{array}
\right)
\end{equation*}
is equal to zero. Since previous determinant is equal to $b_{11}$, we obtain $b_{11}=0$.

\smallskip\noindent
Next, let us consider the following $3$-parameter family of
Lagrangian planes
\begin{equation*}
L=\langle  e_1,e_2+p_{22}e^2 + p_{23}e^3, e_3+p_{23}e^2 + p_{33}e^3, e^4, e^5, \dots , e^n    \rangle,
\end{equation*}
where $p_{22}$, $p_{23}$ and $p_{33}$ are arbitrary real constants.

\smallskip\noindent
Each of such Lagrangian plane intersects $D$, which implies that
the determinant of the following matrix
\begin{equation*}
\left(
\begin{array}{ccccc|cccccc}
1 & 0 & 0 & \cdots  & 0 & 0 & 0 & 0 & 0 & \cdots  & 0 \\
0 & 1 & 0 & \cdots  & 0 & 0 & p_{22} & p_{23} & 0 & \cdots  & 0 \\
0 & 0 & 1 & \cdots  & 0 & 0 & p_{23} & p_{33} & 0 & \cdots  & 0 \\
0 & 0 & 0 & \cdots  & 0 & 0 & 0 & 0 & 1 & \cdots  & 0 \\
\vdots  & \vdots  & \vdots  & \ddots  & \vdots  & \vdots  & \vdots  & \vdots
& \vdots  & \ddots  & \vdots  \\
0 & 0 & 0 & \cdots  & 0 & 0 & 0 & 0 & 0 & \cdots  & 1 \\
\hline
1 & 0 & 0 & \cdots  & 0 & b_{11}=0 & b_{12} & b_{13} & b_{14} & \cdots  &
b_{1n} \\
0 & 1 & 0 & \cdots  & 0 & b_{21} & b_{22} & b_{23} & b_{24} & \cdots  &
b_{2n} \\
\vdots  & \vdots  & \ddots  & \ddots  & \vdots  & \vdots  & \vdots  & \vdots
& \vdots  & \vdots  & \vdots  \\
0 & 0 & 0 & \ddots  & 0 & b_{n-1\;1} & b_{n-1\;2} & b_{n-1\;3} & b_{n-1\;4}
& \cdots  & b_{n-1\;n} \\
0 & 0 & 0 & \cdots  & 1 & b_{n1} & b_{n2} & b_{n3} & b_{n4} & \cdots  &
b_{nn}
\end{array}
\right)
\end{equation*}
vanishes for each choice of $p_{22}$, $p_{23}$, $p_{33}$. But the previous determinant is equal to
\begin{equation}\label{eq.det.ric}
\det\left(\begin{array}{ccc}
0 & b_{12} & b_{13} \\
b_{21} & b_{22}-p_{22} & b_{23}-p_{23} \\
b_{31} & b_{32}-p_{23} & b_{33}-p_{33}
\end{array}
\right)
\end{equation}
So, if we choose $p_{22}=b_{22}$, $p_{23}=b_{23}$ and
$p_{33}=b_{33}$, we get the following equation
\begin{equation}
b_{21}b_{13}(b_{32}-b_{23})=0.
\end{equation}
\begin{itemize}
\item[\emph{First case}:] $b_{21}=0$.

\smallskip\noindent
In this case \eqref{eq.det.ric} is equal to
$b_{31}\big(b_{12}(b_{23}-p_{23})-b_{13}(b_{22}-p_{22})\big)$. If
$b_{31}=0$ we obtain $b_{21}=b_{31}=0$.

\smallskip\noindent
If $b_{31}\neq 0$, then $\big(b_{12}(b_{23}-p_{23})-b_{13}(b_{22}-p_{22})\big)=0$ for any $p_{22}$, $p_{23}$, which implies $b_{12}=b_{13}=0$.

\item[\emph{Second case}:] $b_{21}\neq 0$, $b_{13}=0$.

\smallskip\noindent
This case is analogous to the first case, and then we shall not discuss it.

\item[\emph{Third case}:] $b_{21}\neq 0$, $b_{13}\neq 0$, $b_{23}=b_{32}$.

\smallskip\noindent
If we put in matrix \eqref{eq.det.ric} $p_{22}=b_{22}-1$, $p_{23}=b_{23}$, $p_{33}=b_{33}$, then its determinant is equal to $b_{13}b_{31}$. Since this determinant vanishes,we obtain $b_{31}=0$, i.e. in the same situation of first case.
\end{itemize}
So, we arrived to the following alternative (that we call $\beta_{23}$):
\begin{flushleft}
$(\beta_{23}):\qquad \qquad  b_{12}=b_{13}=0\,,\quad \text{or } \quad b_{21}=b_{31}=0$.
\end{flushleft}
In addition, the above reasoning for indices $2,3$, can be repeated for any couple $i,j=2\dots n$. In this way, for any $i,j$,
\begin{flushleft}
$(\beta_{ij}):\qquad \qquad  b_{1i}=b_{1j}=0\,,\quad \text{or } \quad b_{i1}=b_{j1}=0$.
\end{flushleft}
The collection of alternatives $(\beta_{ij})$ implies
\begin{equation*}
\begin{array}{c}
\text{(A)}\,\,\,\,\,b_{12}=b_{13}=b_{14}=\dots = b_{1n}=0 \\
\text{or}\\
\text{(B)}\,\,\,\,\,b_{21}=b_{31}=b_{41}=\dots =b_{n1}=0
\end{array}
\end{equation*}
Indeed if, for example, $b_{21}\neq 0$, then $(\beta_{1,j})$ implies $b_{12}=b_{1j}=0$, $j=3\dots n$. In other words (A) holds.

\smallskip\noindent
By taking into account that $b_{11}=0$, (A) means $e\in D$ and (B) means $e\in D^\perp$.

\smallskip\noindent
The condition is sufficient. We shall prove that $\mathrm{E}_D\subset\mathrm{E}_{D^\bot}$. If $L$ a Lagrangian plane, we have the following equalities:
\begin{equation}\label{seq}
L\cap D^\bot = L^\bot\cap D^\bot = (L\cup D)^\bot = (L+D)^\bot.
\end{equation}
If furthermore $L\in\mathrm{E}_D$, then by definition $L$ non trivially intersects $D$, that implies $\dim (L+D)\leq n-1$. This means that $\dim(L+D)^\bot\geq 1$, and then $L\cap D^\bot\neq 0$. The same argument leads to the proof of the inverse inclusion.
\end{proof}

\begin{corollary}
Up to a factor, there exist only two decomposable $n$-forms $\Omega_D$ and $\Omega_{D^\perp}$ which give the same equation.
\end{corollary}

\begin{remark}\label{rem.dimensione.preservata}
Note that  subspaces   $L \cap D$ and $L \cap D^\bot$ have  the
same dimension.
In fact by \eqref{seq} we have that
\[
\dim(L\cap D^\bot) = \dim(L+D)^\bot = 2n-\dim(L+D)=2n-(n+n-\dim(L\cap D))=\dim(L\cap D).
\]
\end{remark}
As a corollary of Theorem \ref{th.D.equal.D.orth}, we can reconstruct
$D\cup D^\perp$ from the hypersurface $\mathrm{E}_D$.
\begin{corollary}
Let $(V,\omega)$ be a $2n$-dimensional symplectic vector space and  $D\subset V$ be
an $n$-plane. Then
$$
D\cup D^\perp =\{e\,\,|\,\,e^{(1)}\cap D\neq 0\} = \{e\in V\,\,|\,\,e^{(1)}\subset\mathrm{E}_D\}.
$$
\end{corollary}

\subsubsection{Description of the conformal metric $g_{\mathrm{E}_D}$ and of the singular points of ${\mathrm{E}_D}$}\label{sec.metric.on.taut.point}

Below we describe conformal metric $g_{\mathrm{E}_{D}}$. We need the following technical lemma.
\begin{lemma}\label{lemma.rank.adjoint.matrix}
Let $C$ be an $(n\times n)$ matrix and $A$ its classical adjoint matrix. Then
\begin{enumerate}
\item If $C$ is not degenerate
then $A$ is not degenerate;
\item if $\text{rank}(C)<n-1$ then $A=0$;
\item\label{enum.lemma.rank.adj.mat}
if $\text{rank}(C)=n-1$ then $\text{rank}(A)=1$ and $A=||a^ib^j||$,
where $a$ is solution of the equation $C\cdot x=0$ and $b$ is solution
of the equation $C^t\cdot x=0$. In particular if $A^{nn}=a^nb^n=0$ then either
the last column or the last row is zero.
\end{enumerate}
\end{lemma}
\begin{proof}
Let $c_i$ be the rows of matrix $C$ and $a^j$ the columns
of matrix $A$. Then $c_i\cdot a^j=\det(C)\delta_i^j$. This proves $1$.
Claim $2$ is well known. Now we prove claim $3$. From equation
$c_i\cdot a^j=\det(C)\delta_i^j$ it follows that vectors $a^j$ are solutions
to equation $C\cdot x=0$ and then they are proportional to some solution $a$.
Changing columns and rows in matrices $C$ and $A$, we prove that vector
$b=(b^1,b^2,\dots , b^n)$ is a solution to equation $C^t\cdot x=0$.
\end{proof}
\begin{proposition}\label{prop.metric.ED.point}
Let $\mathrm{E}_D$ be the hypersurface of $\mathcal{L}(V)$  associated with $n$-plane
\eqref{eq.appoggio.D} and $L=L_P=\langle w_i=e_i+p_{ij}e^j \rangle \in  \mathrm{E}_D$. Then the conformal  metric
$g_{\mathrm{E}_{D}}$ in $L^*$ is given by
$$
g_{\mathrm{E}_{D}} = A^{ij}\,w_i\vee w_j
$$
where $A=||A^{ij}||$ is the classical adjoint matrix of matrix
$(P-B)$.\\
Moreover
\begin{enumerate}
\item $A=0$ if $\text{rank}\,(P-B)<n-1$;
\item $A=||a^ib^j||$ if $\text{rank}\,(P-B)=n-1$ where $(P-B)\cdot a=0$
and $(P-B^t)\cdot b=0$. In particular
 \begin{enumerate}
 \item $g_{\mathrm{E}_D}=a\vee b$, $a=a^iw_i$, $b=b^iw_i$;
 \item matrix $\frac{1}{2}(A+A^t)$
of the symmetric form $g_{\mathrm{E}_{D}}$ has rank equal to $1$
if $B=B^t$ and rank equal to $2$ if $B\neq B^t$.
 \end{enumerate}
\end{enumerate}
\end{proposition}
\begin{proof}
Since
\[
\frac{\partial}{\partial p_{ij}}\big(\det(P-B)\big)=\left\{
\begin{array}
[c]{cc}%
A^{ii} & \text{if }i=j\\
A^{ij}+A^{ji} & \text{if }i\neq j
\end{array}
\right.
\]
then
$$
g_{\mathrm{E}_{D}}(\eta,\eta) = \sum_{i\leq j}
\frac{\partial}{\partial p_{ij}}\big(\det(P-B)\big) \eta_i\eta_j =
\sum_{i,j} A^{ij}\eta_i\eta_j = \frac{1}{2}\sum
(A^{ij}+A^{ji})\eta_i\eta_j.
$$
So the matrix of symmetric bilinear form is the symmetrization of
the matrix $A$. This proves the first part of proposition.

\smallskip\noindent
The second part follows from Lemma \ref{lemma.rank.adjoint.matrix}.
\end{proof}

\begin{definition}
A point $L\in\mathrm{E}_D$ is called \textbf{singular} if $\dim(L\cap D)\geq 2$ and \textbf{regular} otherwise. The set of regular points of $\mathrm{E}_D$ will be denotes by $\mathrm{E}^{\mathrm{reg}}_{D}$.
\end{definition}
Now we give a criterion to distinguish singular points.
\begin{proposition}\label{prop.singular.point}
A point $L_P\in\mathrm{E}_D$ is singular iff the differential of
$\det(P-B)$ at $L$ vanishes, that is if the metric
 $g_{\mathrm{E}_D}$ vanishes at $L$.
\end{proposition}
\begin{proof}
We have that
\[
\dim(L\cap D)=k~\Longleftrightarrow~\text{rank}\,(P-B)=n-k,
\]
where $L\in\mathrm{E}_D$. If $k\geq 2$, then $\text{rank}\,
(P-B)\leq n-2$, which implies that its adjoint matrix vanishes in
view of Lemma \ref{lemma.rank.adjoint.matrix}. Then
$\frac{\partial}{\partial p_{ij}}\big(\det(P-B)\big)=0$ at the
point $L$ and $g_{\mathrm{E}_D}|_L=0$ (see also the proof of
Proposition \ref{prop.metric.ED.point}).
\end{proof}

\smallskip\noindent
In view of the definition of singular points, taking into account Remark \ref{rem.dimensione.preservata}, Theorem \ref{th.D.equal.D.orth} restricts to regular points, more precisely we have the following results.
\begin{corollary}
Let $(V,\omega)$ be a $2n$-dimensional symplectic vector space. Let $D$ and $\widetilde{D}$ be $n$-dimensional planes of $V$. Then
$$
\mathrm{E}^{\mathrm{reg}}_{\widetilde{D}}=\mathrm{E}^{\mathrm{reg}}_{D}\,\,\Longleftrightarrow \widetilde{D}=D \text{\,\,\,or\,\,\,} \widetilde{D}=D^\perp.
$$
\end{corollary}

\subsubsection{Description of $\mathrm{E}_{D}$ in terms of its characteristics} \label{sec.characterization}

The theorem below describes characteristic $(n-1)$-dimensional subspaces for hypersurfaces of type $\mathrm{E}_{D}$.
\begin{theorem}\label{th.app}
Let $D$ and $\Omega_D$ be as in Proposition \ref{prop.n.forms}.
Let also $H\subset V$ be an $(n-1)$-dimensional isotropic subspace
and $H^{(1)} = \{ L_t \}$. Then the following conditions are
equivalent:
\begin{enumerate}
\item $H\subset L_0 $ is characteristic for $\mathrm{E}_D$ at $L_0 \in \mathrm{E}_D$;
\item $H^{(1)}\subset\mathrm{E}_D$;
\item $\textrm{vol}_t\,\lrcorner\,\Omega_D=0$, where $\textrm{vol}_t$ is a volume element of $L_t$;
\item $L_t\cap D\neq 0$ for all $t$;
\item $H$ has non trivial intersection with $D$ or $D^\perp$.
\end{enumerate}
\end{theorem}
\begin{proof}
Equivalence $1\Leftrightarrow 2$ is Theorem
\ref{th.strongly.char}, taking into account that
$\mathrm{E}_D=\mathrm{E}_{\Omega_D}$.

\smallskip\noindent
Properties $3$ and $4$ are by definition an alternative ways to write property $2$.

\smallskip\noindent
Now we prove equivalence $2\Leftrightarrow 5$. Let $H$ be characteristic for $\mathrm{E}_D$ at $L$, (so, is also strongly characteristic and then any Lagrangian plane which contains $H$, intersects non trivially $D$). We want derive that $H$ has non trivial intersection with $D$ or $D^\perp$.

\noindent
Let us assume that $H\cap D=0$; we will show that $H\cap
D^\perp\ne 0$.

\noindent
We can take a symplectic basis $\{e_i,e^i\}$ such that
$H=\langle e_1,\dots,e_{n-1}\rangle$ and $L=\langle e_1,\dots,e_{n-1},e_n\rangle$. By
hypothesis, $L\cap D\ne 0$, so that the unique possibility is that $L\cap D$ is generated by a vector $e_n+\sum_{i=1}^{n-1}\alpha_i
e_i$. By a change of the basis we can suppose that this generator
is $e_n$ (in particular, $e_n\in D$). Now, the Lagrangian planes $L_t:=\langle e_1,\dots, e_{n-1}, e_n+te^n\rangle$
have non trivial intersections with $D$. Indeed, by the same
reasoning as above, the intersection $L_t\cap D$, $t\ne 0$, must
be generated by a vector of the form
$$
e_n+te^n+\sum_{i=1}^{n-1}\alpha_i(t)e_i=e_n+t
\left( e^n+\sum_{i=1}^{n-1}\frac{\alpha_i(t)}t e_i\right)
$$
Taking into account that $e_n\in D$, we have
$$
e^n+\sum_{i=1}^{n-1}\frac{\alpha_i(t)}t\, e_i\in D,
$$
If we take two different values $t,\overline t$ we have that
$$
\sum_{i=1}^{n-1}\left(\frac{\alpha_i(t)}t-\frac{\alpha_i(\overline t)}{\overline t}\right)  e_i\in D\cap H=0,
$$
so that
$$
v_n:=e^n+\sum_{i=1}^{n-1}\frac{\alpha_i(t)}t\, e_i
$$
does not depend on $t$. A new change of coordinates allow us to take $e^n=v_n$ so that,
$$
L_t\cap D=\langle e_n+te^n\rangle;
$$
in particular, $D\supset \langle e_n,e^n\rangle$ and $D^\perp\subset \langle e_n,e^n\rangle^\perp$. Also,
$H\subset\langle e_n,e^n\rangle^\perp$ and a computation of dimensions gives us
$$
\dim D^\perp\cap H=\dim D^\perp+\dim H-\dim (D^\perp+H)\ge
    n+(n-1)-(2n-2)=1,
$$
because $D^\perp+H\subset \langle e_n,e^n\rangle^\perp$. Finally, $H\cap D^\perp\ne 0$, as we wanted.
\end{proof}

\begin{remark}
Claims  $1$, $2$, $3$ of the  theorem remain equivalent also for a
hypersurface  $\mathrm{E}_\Omega$, associated with   any $n$-form
$\Omega\in\Lambda^n(V^*)$.
\end{remark}
Bringing together Theorems \ref{th.rotating.debole}, \ref{th.D.equal.D.orth}, \ref{th.app} and Proposition \ref{prop.metric.ED.point}, in the theorem below we will summarize the main results regarding the hypersurfaces of type $\mathrm{E}_{D}$ by putting in evidence how to describe them in terms of their characteristics.
\begin{theorem}\label{theorem.main.2}
Let $\mathrm{E}^{\mathrm{reg}}_D$ be the set of regular point of $\mathrm{E}_D$. Then
\begin{itemize}
\item A hyperplane $H$ of $L\in\mathrm{E}^{\mathrm{reg}}_D$ is characteristic for $\mathrm{E}^{\mathrm{reg}}_D$ at $L$
iff it contains one
of the following straight lines:
$$
\ell_L:=L\cap D \quad \text{or} \quad \ell'_L:=L\cap D^\perp.
$$
Then, if $\ell_L\neq\ell'_L$, there are two $(n-2)$-parametric
families $H(t_1,\dots,t_{n-2})$
and $H'(t_1,\dots,t_{n-2})$ of characteristic hyperplanes in $L$: one contains
$$
\ell_L = \bigcap_{t_1,\dots,t_{n-2}}H(t_1,\dots,t_{n-2})
$$
and another contains
$$
\ell'_L=\bigcap_{t_1,\dots,t_{n-2}}H'(t_1,\dots,t_{n-2}).
$$
If $\ell_L=\ell_L'$ then these two families coincide.
\item The conformal metric of $\mathrm{E}^{\mathrm{reg}}_D$ is decomposable and is given by
$$
\left(g_{\mathrm{E}^{\mathrm{reg}}_D}\right)_L=\ell_L\vee\ell_L'.
$$
\item For any line $\ell\subset D$ there exists   $L\in \mathrm{E}^{\mathrm{reg}}_D$
such that $\ell=\ell_L=L\cap D$. Hence
$$
D=\bigcup_{L\in\mathrm{E}_D} \ell_L\,,\quad D^\perp=\bigcup_{L\in\mathrm{E}_D}\ell'_L.
$$
\end{itemize}
\end{theorem}

\section{Contact manifolds and scalar PDEs of $1^{st}$ order}\label{sec.ContactManifolds}

\begin{definition}\label{completa_non_integrablita}
A $(2n+1)$-dimensional smooth manifold $M$ endowed with a completely
non-integrable  codimension one distribution  $\mathcal{C}$ is
called a \textbf{contact manifold}. A diffeomorphism $\Psi$ of $M$ which preserves $\mathcal{C}$ is called a \textbf{contact transformation}.
\end{definition}
Locally $\mathcal{C} = \Ker\, \theta$, where the contact form
$\theta$ is defined up to a conformal factor. There exist coordinates $(x^i, z, p_i)$, $i=1,\dots, n$ such that \begin{equation}\label{coordenadas de contacto}
\theta = dz -  p_i dx^i.
\end{equation}
Such  coordinates  are called \textbf{contact (or Darboux) coordinates}.
Locally defined vector fields
\begin{equation}\label{eq.basis.of.C}
\widehat{\partial}_{x^{i}}\overset{\text{def}}{= }
\partial_{x^{i}}+p_{i}\partial_{z},\quad\partial_{p_{i}},\qquad
i=1,\dots,n.
\end{equation}
span the contact distribution $\mathcal{C}$. We remark that, in view of the complete non-integrability of $\mathcal{C}$, the contact form $\theta$ cannot depend on $k$ $1$-forms, with $k\leq n$. From now on, for
simplicity, we will assume that the contact form $\theta$ is
globally defined. The 2-form $ d \theta $ is non degenerate on
$\mathcal{C}_m$, $\forall\,m\in M$. We will consider the symplectic structure
$$
\omega = d\theta|_{\mathcal{C}}
$$
in the distribution $\mathcal{C}$. A contact transformation induces a
conformal transformation both of $\theta$ and of $\omega$, so that with any contact manifold a conformal symplectic structure on the contact distribution is associated.

\smallskip\noindent
Recall that a \textbf{Legendre transformation} is a local
contact transformation $(x^i,z,p_i) \to (x'^i,z',p'_i)$ defined  by
\begin{equation*}
x'^i=p_{i},\,\, z'=z-p_{i}x^{i},\,\, p'_{i}=-x^{i},\,\,
i=1,\dots,n.
\end{equation*}
The action of such transformation on vector fields interchanges the
roles of $\widehat{\partial}_{x^{i}}$ and $\partial_{p_{i}}$; indeed,
\begin{equation}\label{legendre.1}
 \partial_{z}  \mapsto \partial_{z'}%
,\quad  \widehat{\partial}_{x^{i}}  \mapsto
-\partial_{p'_{i}},\quad  \partial_{p_{i}}  \mapsto \widehat{\partial}_{x'^{i}}.
\end{equation}
Sometimes it is useful to define a ``partial'' Legendre transformation. For
instance, we can divide the indices $i=1,\dots,n$ into $\alpha=1,\dots,m$ and
$\beta=m+1,\dots,n$ and define
\begin{equation}\label{legendre.2}
z'=z-p_{\alpha}x^{\alpha},\,\,
x'^{\alpha}=p_{\alpha},\,\,
p'_{\alpha}=-x^{\alpha},\,\,
x'^{\beta}=x^{\beta},\,\,
p'_{\beta}=p_{\beta},\,\,\, \alpha=1,\dots,m\,,\,\,\,\beta=m+1,\dots,n
\end{equation}
which also defines a contact transformation. In this case, only the first $m$
coordinates $x^{\alpha}$ and $p_{\alpha}$ are interchanged (joint the
corresponding partial derivatives).

\subsection{Cartan and Hamiltonian vector fields }

\begin{definition}
Sections $Y \in \Gamma(\mathcal{C})$ are called { \bf Cartan vector
fields}.
\end{definition}
Cartan fields form a $C^\infty(M)$-module and vector fields \eqref{eq.basis.of.C} form a local basis. They do not form a Lie algebra: in fact the formula
\begin{equation*}
(Y\cdot \theta)(X)= d\theta(Y,X)=\omega(Y,X)= \theta([X,Y]),\,\,\,
X,Y,\in \Gamma(\mathcal{C})
\end{equation*}
where we recall that $Y\cdot\theta$ is the Lie derivative of $\theta$ along $Y$, shows that two Cartan fields are orthogonal iff their Lie bracket
is still a Cartan field. It allows to  express
$\omega$-orthogonality in $\mathcal{C}$
 in terms of Lie
derivatives. For example, the orthogonal complement of $Y$ in
$\mathcal{C}$ is described by
\[
Y^{\bot}=\{\theta=0,~Y\cdot\theta=0\}.
\]
In particular, $Y^{\bot}$ is $(2n-1)$-dimensional and contains
$Y$; moreover, any $(2n-1)$-dimensional subdistribution of
$\mathcal{C}$ is of this form. Analogously, if
$\mathcal{D}\subset\mathcal{C}$ is a distribution spanned by
vector fields $Y_{1},\dots,Y_{k}$ then its orthogonal complement
is given by
\[
\mathcal{D}^{\bot}=\{\theta=0,\,Y_{1}\cdot\theta=0,\dots,Y_{k}\cdot\theta=0\}.
\]
The flow generated by a Cartan field $Y$
deforms $\mathcal{C}$, and the sequence of iterated Lie derivatives
\begin{equation}\label{Y_i(theta)}
\theta,\,\,\,\,Y\cdot\theta,\,\,\,\,Y\cdot (Y\cdot\theta),\,\dots,\,\,
\underset{(2n-1)\text{-times}}{\underbrace{Y\cdot (Y\cdot \dots \cdot (Y}}\cdot \theta )\dots )
\end{equation}
gives a measure of this deformation.
\begin{definition}
The {\bf type} of  a Cartan field $Y$ is defined as the rank of system
(\ref{Y_i(theta)}).
\end{definition}
Let us fix a contact form $\theta$; the Reeb vector field $Z$ is defined by conditions
$$
\theta(Z) =1,\,\,\, Z\,\lrcorner\,\omega =0.
$$
It depends on the choice of $\theta$. We denote by  $Z^0 \subset
\Lambda^1(M) $ the annihilator of $Z$ in the space of 1-forms. In a
contact chart (\ref{coordenadas de contacto}), $Z=\partial
/\partial z$ and the following decomposition holds:
\begin{equation*}
TM\simeq\langle
Z\rangle\oplus\mathcal{C},\qquad v\mapsto \theta(v)Z+
(v-\theta(v)Z)
\end{equation*}
or, dually,
\begin{equation}\label{suma directa 2}
T^{*}M\simeq\langle \theta\rangle\oplus
Z^0,\qquad \alpha\mapsto\alpha(Z)\theta+(\alpha-\alpha(Z)\theta).
\end{equation}
The map
\begin{equation*}
\chi:\Gamma(\mathcal{C})\to   Z^0, \quad   Y \mapsto Y \cdot
\theta = Y\,\lrcorner\,d\theta
\end{equation*}
is an isomorphism of $C^{\infty}(M)$-modules. So any 1-form
$\alpha\in\Lambda^{1}(M)$  defines  a Cartan vector field
\[
Y_{\alpha}\overset{\text{def}}{=}\chi^{-1}\left(
\alpha-\alpha(Z)\theta\right)
\]
(see the direct sum (\ref{suma directa 2})). In other words, $Y_{\alpha}%
\in\mathcal{C}$ is determined by the relation%
\begin{equation*}
Y_{\alpha}\cdot\theta=Y_{\alpha}\lrcorner d\theta=\alpha-\alpha(Z)\theta.
\end{equation*}
So any  Cartan vector  field has  the form $Y_{\alpha}$  and $1$-form $\alpha$ is canonically defined up to
adding  a  form proportional to $ \theta$.  We have
\begin{equation*}
 Y_\alpha \cdot \theta = \alpha - \alpha(Z) \theta,\,\,\,  Y_{f \alpha} = f Y_{\alpha}, \,\, \forall f \in C^{\infty}(M).
\end{equation*}
If we choose a different generator $\theta^{\prime}=\lambda \theta$ we have that
\[
X^{\prime}_{\alpha}=\frac1\lambda X_{\alpha};
\]
in particular, although $X_{\alpha}$ depends on the choice of $\theta$, its direction does not change.

\begin{definition}
A vector field $Y_f:=Y_{df}$ is called a {\bf Hamiltonian vector field}.
\end{definition}
In contact coordinates $(x^{i},z,p_{i})$ a Hamiltonian vector field can be written as
\[
Y_{f}=\sum_{i=1}^{n} \,\, \partial_{p_{i}}(f)\widehat{\partial}_{x^{i}} -
\widehat{\partial}_{x^{i}}(f)\partial_{p_{i}}.
\]
In particular
$
Y_{x^{i}}=-\partial_{p_{i}}, \,\,\, Y_{z}=-\sum_{i=1}^{n} p_{i}\partial_{p_{i}%
},\,\,\, Y_{p_{i}}=\widehat{\partial}_{x^{i}}.
$

\medskip\noindent
From the above definition, next lemma easily follows.
\begin{lemma}
A Hamiltonian vector field $Y_f$ satisfies the following equalities
\begin{equation}\label{eq.proprieta.Yf}
df(Y_f)=Y_f(f)=0\,, \quad \theta(Y_f)=0\,, \quad Y_f\cdot\theta=df-\frac{\partial f}{\partial z}\theta.
\end{equation}
\end{lemma}

\begin{remark}\label{rem.Yf.char.symm}
Previous lemma implies that $Y_{f}$ is a characteristic symmetry for the distribution
$Y_{f}^{\bot}=\{\theta=0,\,df=0\}$. In other words, $Y_{f}$ coincides with the
classical characteristic vector field of the first order equation $f(x^i,z,p_i)=0$ where $p_i = \partial z/\partial x^i$. Also, properties \eqref{eq.proprieta.Yf} easily imply that $Y_{f}$ is a vector field of type $2$.
\end{remark}
\begin{definition}
Two functions $f$ and $g$ on $M$ are in \textbf{involution} if $\omega(Y_f,Y_g)=0$ (or equivalently, if $Y_f(g)=0$).
\end{definition}
\begin{lemma}\label{lem.involution.non.so}
Two functions $f$ and $g$ on $M$ are in involution iff the distribution $\langle Y_f,Y_g \rangle$ is integrable.
\end{lemma}
\begin{proof}
We have that
\begin{equation}\label{eq.f.g.involution}
\omega(Y_f,Y_g)=-\theta([Y_f,Y_g]).
\end{equation}
Now let us suppose that $f$ and $g$ are in involution. Then previous equality implies that $[Y_f,Y_g]\in\Gamma(\mathcal{C})$. On the other hand it is easy to see that
$$
[Y_f,Y_g]\cdot\theta=\lambda Y_f\cdot\theta + \mu Y_g\cdot\theta + \nu\theta
$$
for some functions $\lambda$, $\mu$, $\nu$. In this way
$$
([Y_f,Y_g]-\lambda Y_f - \mu Y_g)\cdot\theta=\nu\theta,
$$
which implies that
$$
[Y_f,Y_g]-\lambda Y_f - \mu Y_g=0,
$$
since a non-trivial Cartan field cannot be an infinitesimal symmetry of $\mathcal{C}$.

\smallskip\noindent
If $\langle Y_f,Y_g \rangle$ is integrable, then equality \eqref{eq.f.g.involution} implies that $f$ and $g$ are in involution.
\end{proof}

\smallskip\noindent
The theorem below is extracted from \cite{Munoz}.
\begin{theorem}\label{th.munoz}
Any set $(f_1,\dots,f_k)$ of $k$ functions on the contact manifold $M$ which are in involution can be extended to a contact chart.
\end{theorem}
\begin{proof}
By Lemma \ref{lem.involution.non.so}, distribution $\mathcal{P}=\langle Y_{f_1}\,\dots,Y_{f_k}\rangle$ is integrable. In particular $\mathcal{P}$ is isotropic and $k\leq n$. If $k<n$, in view of Lemma \ref{lem.involution.non.so}, we can take a first integral $f_{k+1}$ of $\mathcal{P}$ such that distribution $\langle Y_{f_1}\,\dots,Y_{f_{k+1}}\rangle$ is $(k+1)$-dimensional and integrable. By iterating this process, we get an $n$-dimensional integrable distribution $\langle Y_{f_1}\,\dots,Y_{f_n}\rangle=\langle d{f_1}=\cdots =d{f_n}=\theta=0\rangle$. So there exists a function $f_0$ such that $\theta=\sum_{i=0}^n a_i df_i$. Then
$$
z=f_0\,,\,\,\, x^i=f_i\,,\,\,\, p_i=-\frac{a_i}{a_0}\,,\,\,\, i=1,\dots,n
$$
gives a contact chart on $M$.
\end{proof}

\subsection{Integral submanifolds of the contact distribution}\label{sec.int.subman.contact}

Recall that an integrable subdistribution of $\mathcal{C}$  is $\omega$-isotropic, hence it has dimension $\leq n$.
As is well know, any $n$-dimensional integral distribution of $\mathcal{C}$, if parametrizable by $(x^1,\dots,x^n)$, is of the form:
$$
z=g(x^1,\dots,x^n),\,\quad p_i=\frac{\partial g}{\partial x^i}(x^1,\dots,x^n).
$$
Integral distributions of $\mathcal{C}$ of dimension $(n-1)$ are described below. The following lemma is a version of classical method of characteristics.
\begin{lemma}\label{lemma.char.flow}
Let $N$ be an integral submanifold of $\mathcal{C}$, $f\in C^\infty(M)$ such that $f|_N=0$ and $\varphi_t$ be the local flow of $Y_f$. Then $\bigcup_t\varphi_t(N)$ is a solution of $\theta=0$ and also of $f=0$.
\end{lemma}
\begin{proof}
In view of Remark \ref{rem.Yf.char.symm}, the local flow $\varphi_t$ of $Y_f$ preserves solutions of the Pfaff system $\{\theta\,,df\}$. So $\bigcup_t\varphi_t(N)$ is a solution of both $\theta=0$ and $f=0$.
\end{proof}
\begin{proposition}
An $(n-1)$-dimensional submanifold $N$ is an integral submanifold of $\mathcal{C}$ iff it is a hypersurface of an $n$-dimensional integral submanifold (of $\mathcal{C}$).
\end{proposition}
\begin{proof}
Of course the condition is sufficient. We prove that it is also necessary. Let us consider a function $f$ on $M$ such that $f|_N=0$ and $({Y_f})_m\cap T_mN=0$ for any $m\in N$. Such a function always exists. In fact, if
$$
N=\{f_1=0,\dots,f_{n+2}=0\}
$$
then the $(n+2)$ Hamiltonian vector fields $Y_{f_i}$ cannot be simultaneously tangent to $N$ for dimensional reasons. The proposition follows in view of above lemma.
\end{proof}
\begin{corollary}\label{cor.caratterizzazione.cauchy.datum}
Let $N$ be an integral $(n-1)$-dimensional submanifold of
$\mathcal{C}$. Then for any point of $N$ there exists a neighborhood in $N$
which is described by
$$
\left\{x^1,x^2,\dots,x^n=0,\,z=\phi(x^1,\dots,x^{n-1}),\,p_1=\frac{\partial\phi}{\partial x^1},\dots,p_{n-1}=\frac{\partial\phi}{\partial x^{n-1}},p_n=\phi_n(x^1,\dots,x^{n-1})\right\}
$$
w.r.t. some local contact coordinates $(x^i,z,p_i)$ of $M$ for certain functions $\phi$ and $\phi_n$. Furthermore, we can select a new contact chart $(\overline{x}^i,\overline{z},\overline{p}_i)$ by taking $\overline{z}=z-\phi$ so that in this new chart $N$ is described by
$$
\left\{x^1,x^2,\dots,x^n=0,\,z=0,\,p_1=0,\dots,p_{n-1}=0,p_n=\phi_n(x^1,\dots,x^{n-1})\right\}.
$$
\end{corollary}

\subsection{Scalar PDEs of $1^{st}$ order and methods of characteristics}

\begin{definition}
A scalar first order partial differential equation ($1^{st}$ order PDE) with one unknown function and
$n$ independent variables is a hypersurface $\mathcal{F}$ of a $(2n+1)$%
-dimensional contact manifold $(M,\mathcal{C})$. A solution of $\mathcal{F}$
is, by definition, an integral manifold of $\mathcal{C}$ contained in
$\mathcal{F}$.
\end{definition}
Clearly the dimension of a solution of $\mathcal{F}$ is less or equal to $n$, as it is also an integral
manifold of $\mathcal{C}$.
In terms of coordinates, $\mathcal{F}$  can be described as a zero
level set
$$
M_f:=\{f(x^{i},z,p_{i})=0\}
$$
of a function $f$. A solution $\Sigma$  parametrized by
$x^{1},\dots,x^{n}$ can be written as
\begin{equation*}
\Sigma\equiv%
\begin{cases}
\displaystyle{z=\phi(x^{1},\dots,x^{n})}\\
\\
\displaystyle{p_{i}=\frac{\partial\phi}{\partial x^{i}}(x^{1},\dots,x^{n})}%
\end{cases}
\end{equation*}
where   the  function $\phi$ satisfies
\[
f\left( x^{i},\phi,\frac{\partial\phi}{\partial x^{i}}\right) =0,
\]
which coincides with the classical notion of solution.

\begin{remark}
The role of coordinates ``$x^{i}$'' as independent variables is purely
external. A contact transformation can change the aforesaid role. For
instance, a total or partial Legendre transformation (see \eqref{legendre.1} and \eqref{legendre.2}) can be used in order to consider ``$p_{i}%
$'' coordinates (all or some of them) as new independent variables.
\end{remark}

\begin{definition}
A \textbf{Cauchy datum} for a first order PDE $M_f=\{f=0\}$, $f\in C^\infty(M)$, is an $(n-1)$-dimensional
 integral submanifold of $\mathcal{C}$ included in $M_f$. It is called \textbf{non-characteristic} if
 it is transversal to the Hamiltonian vector field $Y_f$.
\end{definition}
\begin{remark}
The name ``non-characteristic'' is justified since $Y_f$ coincides with the classical characteristic vector field of first order PDE $M_f$ (see Remark \ref{rem.Yf.char.symm}). The name ``Cauchy datum" is justified in view of the following fact: in the case that $M$ is the space $J^1(\mathbb{R}^n)$ of $1$-jets of functions on $\mathbb{R}^n$, an $(n-1)$-dimensional submanifold $N'$ of $\mathbb{R}^n$ can be prolonged in a unique way to a Cauchy datum $N$ for equation $f=0$ without solving any differential equation. In coordinates, if $(x^i,z,p_i)$ is a contact chart on $M=J^1(\mathbb{R}^n)$ and $N'$ is locally described by
$$
N':\,\,\,x^i=\phi^i(t_1,\dots,t_{n-1})\,,\,\,\,z=\phi(t_1,\dots,t_{n-1}),
$$
then
$$
N:\,\,\,x^i=\phi^i(t_1,\dots,t_{n-1})\,,\,\,\,z=\phi(t_1,\dots,t_{n-1}) \,,\,\,\,p_i=\psi_i(t_1,\dots,t_{n-1}),
$$
where functions $\psi_i$ are uniquely determined by the system of $n$ algebraic equations
$$
\left\{
\begin{array}{l}
0=(dz-p_idx^i)|_N=\left(\frac{\partial\phi}{\partial t_h} - \psi_i(t)\frac{\partial\phi^i}{\partial t_h}\right)\,dt_h\\
\\
0=f|_N=f(\phi^i(t),\phi(t),\psi_i(t))
\end{array}
\right.
$$
\end{remark}
Now, let us consider a given Cauchy datum $N$ for the equation $M_f=\{f=0\}$. Then, by Lemma \ref{lemma.char.flow}, manifold $\Sigma=\bigcup_t \varphi_t(N)$, where $\varphi_t$ is the local flow of the Hamiltonian vector field $Y_f$, is a solution of $f=0$. This solution is, locally, the unique which contains $N$, because by Lemma \ref{lemma.char.flow} $Y_f$ is tangent to any maximal solution of $M_f$. In more concrete terms, construction of solutions of first order PDE $f=0$ goes along the following steps:
\begin{enumerate}
\item take a non-characteristic Cauchy datum $N$;
\item integrate vector field $Y_f$;
\item take the set $\Sigma$ of integral curves of $Y_f$ crossing $N$.
\end{enumerate}
The above method is called the \textbf{method of characteristics} (see also \cite{Vinogradov and Co}).

\section{Characteristics of general $2^{nd}$ order PDEs, general MAEs and MAEs of Goursat type}\label{subsec.PDE}

\subsection{Prolongation of a contact manifold and  its  submanifolds}\label{subsec.PDE.parte.prima}

Let $(M,\mathcal{C})$ be a contact manifold. We recall that it defines a conformal symplectic structure $\omega=d\theta|_{\mathcal{C}}$ on $\mathcal{C}$, where $\theta$ is any $1$-form such that $\Ker(\theta)=\mathcal{C}$. We also recall that
$\mathcal{L}(\mathcal{C}_m)$ denotes  the Lagrangian Grassmannian of $(\mathcal{C}_m, \omega_m)$, $m\in M$.
\begin{definition}
The prolongation of a contact manifold $(M,\mathcal{C})$ is the fiber bundle $\pi: M^{(1)}\to M$ where
$$
M^{(1)}=\bigcup_{m\in M}\mathcal{L}(\mathcal{C}_m)
$$
is the set of all Lagrangian planes of the contact distribution.
\end{definition}
Points of $M^{(1)}$ are Lagrangian planes of $(\mathcal{C}_m,\omega_m)$, $m\in M$: a generic point of
$M^{(1)}$ will be denoted either by $m^1$ or by $L_{m^1}$ so that the tautological bundle
$$
{\mathcal{T}}(M^{(1)})=\{(m^1,v)\,\,|\,\,v\in L_{m^1}\} \to
M^{(1)},\,\,  (m^1,v) \mapsto m^1
$$
over $M^{(1)}$ is  well defined.

\smallskip\noindent
Obviously all that we said in Sections \ref{sec.Lagr.Grass.of.sympl.sp} and \ref{sec.sub.of.Lagr.Grass}
can be applied to the fibers of $M^{(1)}$, i.e. to $\mathcal{L}(\mathcal{C}_m)$.

\smallskip\noindent
A system of contact coordinates $(x^i,z,p_i)$ on $M$ induces  coordinates
\begin{equation}
(x^i, z, p_i, p_{ij}= p_{ji} ,\,\,  1 \leq i, j \leq n )
\label{coordinates_prolongation}%
\end{equation}
on $M^{(1)}$ as follows: a point  $m^1 \equiv L_{m^1} \in M^{(1)}$
has coordinates (\ref{coordinates_prolongation}) iff $m=\pi(m^1)=(x^i,x,p_i)$ and the
corresponding Lagrangian plane $L_{m^1}$ is given by:
\begin{equation*}
L_{m^1}=L_P=\langle \widehat{\partial}_{x^{i}}+ p_{ij}\partial_{p_{j}}\rangle\subset\mathcal{C}_{m},
\end{equation*}
where $P=\|p_{ij}\|$, $\widehat{\partial}_{x^{i}}$ are defined in \eqref{eq.basis.of.C} and all vectors are taken in the point $m$. Note that the isotropy condition entails that $p_{ij}=p_{ji}$,
so that the number of  ``second order'' coordinates  $p_{ij}$ is $\frac{n(n+1)}{2}$ and
$\dim M^{(1)}=\frac{1}{2}(n^{2}+5n+2)$.

\smallskip\noindent
An integral submanifold $N$ of the contact manifold
$(M,\mathcal{C})$ (i.e. $TN\subset\mathcal{C}$) is called
\textbf{isotropic}. Note that  $T_mN$ is an isotropic subspace of
$\mathcal{C}_m$, since $\theta|_N=0$ implies
$\omega|_N=d\theta|_N=0$. Maximal ($n$-dimensional) integral
submanifolds of $\mathcal{C}$ are called \textbf{Lagrangian}.

\smallskip\noindent
We define the \textbf{prolongation} $N^{(1)}\subset M^{(1)}$
of a submanifold $N$ of a contact manifold $M$ as the set of
all Lagrangian planes $L$ which are prolongations of the tangent spaces of $T_mN$
 (see \eqref{eq.prolongation.of.U}):
$$
N^{(1)}:=\left\{
\begin{array}{c}
m^1\in M^{(1)}\,\,|\,\, L_{m^1}\supseteq T_mN\cap\mathcal{C}_m,\,\,\text{if}\,\,\dim(N)\leq n\\
\\
m^1\in M^{(1)}\,\,|\,\, L_{m^1}\subseteq T_mN\cap\mathcal{C}_m,\,\,\text{if}\,\,\dim(N)\geq n\\
\end{array}
\right.
$$
If $N$ is an isotropic submanifold, then he natural projection
$\pi_N:N^{(1)}\to N$ is a fibre bundle whose typical fibre is
$U\oplus \mathcal{L}(W)\simeq \mathcal{L}(\mathbb{R}^{2n-2k})$ where
$U$ and $W$ are as in Proposition \ref{prop.prol.U}, with $U=T_mN$
and $V= \mathcal{C}_m$. In particular, if $N$ is a Lagrangian
submanifold, then $N^{(1)}$ consists of tangent spaces of $N$
(which are Lagrangian) and the projection $\pi_N$ is a
diffeomorphism.

\subsection{Characteristic cone and characteristic subspaces of a  PDE $\mathcal{E}$ of $2^{nd}$ order and its conformal metric $g_{\mathcal{E}}$}\label{sec.char}

\begin{definition}
Let $(M,\mathcal{C})$ be a $(2n+1)$-dimensional contact manifold and $M^{(1)}$ its prolongation.
A hypersurface $\mathcal{E}$ of $M^{(1)}$ is called a scalar second order partial differential equation ($2^{nd}$ order PDE) with one unknown function and $n$ independent variables. A solution of $\mathcal{E}$ is a Lagrangian submanifold $\Sigma\subset M$
whose prolongation $\Sigma^{(1)}$ is contained in $\mathcal{E}$.
\end{definition}
As in the first order case, if
$\mathcal{E}=\{F(x^{i},z,p_{i},p_{ij})=0\}$ then a solution
$\Sigma$  parametrized by $x^{1},\dots, x^{n}$, can be written as
\begin{equation*}
\Sigma\equiv%
\begin{cases}
\displaystyle{z=\varphi(x^{1},\dots,x^{n})}\\
\\
\displaystyle{p_{i}=\frac{\partial\varphi}{\partial x^{i}}(x^{1},\dots,x^{n}%
)}\\
\\
\displaystyle{p_{ij}=\frac{\partial^{2}\varphi}{\partial x^{i}\partial x^{j}%
}(x^{1},\dots,x^{n})}%
\end{cases}
\end{equation*}
where  the  function $\varphi$ satisfies the equation
\[
F\left( x^{i},\varphi,\frac{\partial\varphi}{\partial x^{i}},\frac
{\partial^{2}\varphi}{\partial x^{i}\partial x^{j}}\right) =0,
\]
which coincides with the classical notion of solution.

\smallskip\noindent
The restriction of $\pi:M^{(1)}\to M$ to the equation $\mathcal{E}\subset M^{(1)}$ is a fibre bundle whose fibre at $m$ is denoted by $\mathcal{E}_m$:
$$
\mathcal{E}_m:=\mathcal{E}\cap \pi^{-1}(m).
$$
Obviously, all definitions and properties of Section \ref{sec.char.cone.Lagr.Grass} are still valid on
fibres $\mathcal{E}_m$: we can find them just by substituting $\mathrm{E}$ with $\mathcal{E}_m$, $m\in M$.
 Below we resume such properties.

\begin{definition}
A \textbf{Cauchy datum} for a second order PDE is an $(n-1)$-dimensional integral submanifold of the contact
distribution $\mathcal{C}$.
\end{definition}

\begin{definition}
The set
$$
\text{Ch}_{m^1}(\mathcal{E})=
T_{m^1}\mathcal{E}_m\cap T^{1}_{m^1}\mathcal{L}(\mathcal{C}_m)
$$
of rank $1$ (vertical) tangent vectors to the hypersurface $\mathcal{E}$ at $m^1$
is called the \textbf{characteristic cone} of the equation $\mathcal{E}$ at $m^1$. Elements of $\text{Ch}_{m^1}(\mathcal{E})$ are called \textbf{characteristic vectors} for $\mathcal{E}$ at $m^1$. The $1$-dimensional vector space generated by a characteristic vector is called a \textbf{characteristic direction}. A characteristic vector $v$ for $\mathcal{E}$ at $m^1$ is called \textbf{strongly characteristic} if the line $\ell(m^1,v)$ (see the end of Section \ref{sec.Lagr.Grass}) is contained in $\mathcal{E}_m$.
\end{definition}
\begin{definition}
A subspace $U\subset T_mM$ is said to be \textbf{characteristic} for
the equation $\mathcal{E}$ at $m^1$ if $U^{(1)}$ is tangent to
$\mathcal{E}$ at $m^1$. If in addition $U^{(1)}\subset\mathcal{E}$, $U$ is said to be
\textbf{strongly characteristic}. A submanifold $S\subset M$ is said to be characteristic for $\mathcal{E}$
(resp. strongly characteristic) if, for any $m\in S$, $T_mS$ is characteristic at least for a point $m^1\in\mathcal{E}$ (resp. strongly characteristic).
\end{definition}
We would like to underline that previous definitions, in view of Remark \ref{rem.conformal.change}, are invariant under a conformal change of the contact form.

\smallskip\noindent
Remark \ref{rem.rank1.and.hyperplane} explains the relationship between characteristic directions and characteristic subspaces of $\mathcal{E}$. As we did in Section \ref{sec.char.cone.Lagr.Grass}, we can introduce a conformal metric $(g_{\mathcal{E}})_{m^1}=g_{\mathcal{E}_{\pi(m^1)}}$ on $S^2(L^*_{m^1})$ at each point $m^1\equiv L_{m^1}\in M^{(1)}$ and Theorem \ref{th.fab} is still valid mutatis mutandis. In coordinates, a tangent vector to $\mathcal{E}_m$ at $m^1$ having $\dot{P}=||\dot{p}_{ij}||$ as matrix of coordinates is of rank $1$ iff $\dot{p}_{ij}=\eta_i\eta_j$ up to a sign (see also \eqref{eq.rank.1.vectors}). Furthermore, it is characteristics if it satisfies Equation \eqref{eq.characteristics}. A covector $\eta$ is characteristic for $\mathcal{E}$ (see also correspondences \eqref{eq.correspondence.gianni} and \eqref{eq.correspondence.gianni.2}) iff it is isotropic for $g_{\mathcal{E}}$. In view of Theorem \ref{th.rotating.debole}, $(g_\mathcal{E})_{m^1}$ is decomposable iff characteristic hyperplanes of $L_{m^1}$ are divided in two $(n-2)$-parametric families $\mathcal{H}_{m^1}$ and $\mathcal{H}'_{m^1}$ such that
\begin{equation*}
\dim\bigcap_{U\in\mathcal{H}_{m^1}}U=\dim\bigcap_{U\in\mathcal{H}'_{m^1}}U=1
\end{equation*}

\begin{example}
Here we treat the classical case $n=2$. Let $\mathcal{E}=\{F=0\}$ be a second order scalar PDE and $m^1\in\mathcal{E}$ a regular point. Then $\eta=(\eta_1,\eta_2)$ is a characteristic covector for $\mathcal{E}$ at $m^1$ if it satisfies Equation \eqref{eq.characteristics}:
\begin{equation}\label{eq.characteristics.n.2}
\frac{\partial F}{\partial p_{11}}\eta_1^2+\frac{\partial F}{\partial p_{12}}\eta_1\eta_2+\frac{\partial F}{\partial p_{11}}\eta_2^2=0,
\end{equation}
where $\frac{\partial F}{\partial p_{ij}}$ are computed in $m^1$. We note that $(\eta_1^1,\eta_1\eta_2,\eta_2^2)$ is a vector of the characteristic cone of $\mathcal{E}$ at $m^1$.

\noindent
Dually, $v=(v^1,v^2)$ spans a $1$-dimensional characteristic subspace (i.e. a hypersurface of $L_{m^1}$, see correspondence \eqref{eq.correspondence.gianni.2}) for $\mathcal{E}$ at $m^1$ iff
\begin{equation}\label{eq.app2}
\frac{\partial F}{\partial p_{11}}{v^2}^2 - \frac{\partial F}{\partial p_{12}}{v^1v^2} + \frac{\partial F}{\partial p_{22}}{v^1}^2=0
\end{equation}
(compare with \eqref{non_characteristic}). Previous equations have
$2$, $1$ or no real solutions, according to the sign of
\[
\Delta=F_{p_{12}}^{2}-4F_{p_{22}}F_{p_{11}}%
\]
(positive, zero or negative). It follows that left hand side  of
\eqref{eq.characteristics.n.2} and \eqref{eq.app2} are always
decomposable over $\mathbb{C}$. They are decomposable over
$\mathbb{R}$ if $\Delta\geq 0$.
\end{example}

\subsection{Characteristics of general MAEs}\label{sec.char.MAE}

Let $(M,\mathcal{C})$ be a contact manifold and
$\mathcal{I}(\theta)\subset\Lambda^{*}(M)$ be the differential
ideal generated by a contact form $\theta$.
Following V.V. Lychagin (see \cite{KLR,Lychagin contact}), we give the following definition
\begin{definition}\label{def.MA.Lych}
Let $\Omega\in\Lambda^{n}(M) \backslash\mathcal{I}(\theta)$. We associate with $\Omega$ the hypersurface $\mathcal{E}_{\Omega}$ of $M^{(1)}$ defined by
\begin{equation*}
\mathcal{E}_{\Omega}\overset{\textrm{def}}{=}\{{m}^1\in M^{(1)}
\text{ s.t. }\Omega|_{L_{{m}^1}}=0\}=\bigcup_{m\in M}\mathrm{E}_{\Omega_m}\,,
\end{equation*}
where $L_{{m}^1}\subset T_{\pi({m}^1)}M$ is the
Lagrangian plane associated with ${m}^1$ (recall that $\pi$ is the projection of $M^{(1)}$ onto $M$). Equations of this form are called {\bf general Monge-Amp\`{e}re equations}.
\end{definition}
In other words $\mathcal{E}_{\Omega}$ is the differential equation
corresponding to the exterior differential system $\{\theta=0,\,\Omega=0\}$.

\begin{remark}\label{rem.horizontalization}
The correspondence $m^1\in M^{(1)}\mapsto
\Omega|_{L_{m^1}}\in\Lambda^n(L^*_{m^1})$ defines an $n$-form on
the tautological bundle $\mathcal{T}(M^{(1)})$.
\end{remark}

\smallskip\noindent
Two $n$-forms $\Omega, \Omega'$ defines the same equation
$\mathcal{E}_{\Omega}= \mathcal{E}_{\Omega'}$ iff, up to a non vanishing factor, are related by
\begin{equation} \label{forme_equivalenti}
\Omega'=\Omega+\alpha\wedge d\theta+\beta\wedge \theta
 \text{\quad for some $\alpha\in\Lambda^{n-2}(M),\,\,\beta\in\Lambda^{n-1}(M)$}.
\end{equation}

\smallskip\noindent
All results of Section \ref{sec.hyper.Lagr.Grass.assoc.n.forms} can be applied to fibers
${\mathcal{E}_{\Omega}}_m$ just by substituting
$\Omega$ with $\Omega_m$ and ${\mathcal{E}_{\Omega}}_m$ with $\mathrm{E}_{\Omega}$, $m\in M$. In particular, by putting together Theorems \ref{th.strongly.char} and \ref{th.car.hyp}, we obtain the following results.
\begin{theorem}
Let $m^1\in\mathcal{E}_\Omega$. A hyperplane $H \subset L_{m^1}$ is characteristic for the MAE $\mathcal{E}_\Omega$ at $m^1$ if and only if it is strongly characteristic. Moreover,
characteristic hyperplanes are those hyperplanes which are isotropic with respect to  some $n$-form $\Omega'$ equivalent to $\Omega$ in the sense of \eqref{forme_equivalenti}.
\end{theorem}

\subsection{MAEs $\mathcal{E}_{\mathcal D}$ associated with $n$-dimensional subdistributions $\mathcal{D}$ of the contact
distribution and their description in terms of their characteristics}\label{sec.MAE.ass.with.D}

As before, $(M,\mathcal{C})$ is a
$(2n+1)$-dimensional contact manifold and $\theta$ a contact form.
\begin{definition}\label{def.ED}
Let $\mathcal{D}$ be an $n$-dimensional subdistribution of the
contact  distribution  $\mathcal{C}$  of $M$. We associate  with
$\mathcal{D}$ the hypersurface $\mathcal{E}_{\mathcal{D}}$ of $M^{(1)}$ defined by
\[
\mathcal{E}_{\mathcal{D}}\overset{\textrm{def}}{=}\{ m^1\in M^{(1)} \,\,|\,\, L_{m^1}\cap \mathcal{D}_{\pi(m^1)}\neq 0  \}
= \bigcup_{m\in M}\mathrm{E}_{\mathcal{D}_m}.
\]
\end{definition}
\begin{proposition}\label{cor.decomposable}
The equation $\mathcal{E}_{\mathcal{D}}$
defined by an $n$-dimensional subdistribution $\mathcal{D}\subset
\mathcal{C}$ is the MAE associated   with  the n-form
\[
\Omega=\Omega_{\mathcal{D}}:=Y_{1}\cdot\theta\wedge\dots\wedge
Y_{n}\cdot\theta,
\]
where $Y_{i}$ are vector fields generating the orthogonal
distribution $\mathcal{D}^{\bot}$. The converse is also true.
\end{proposition}
\begin{proof}
 Since the subdistribution
$\mathcal{D}\subset \mathcal{C}$ is defined  by the system of
1-forms
\[
\left\{
\begin{array}{r}
\theta=0 \\
Y_i\cdot\theta=0
\end{array}
\right.
\]
where vector fields $Y_{i}$ generate $\mathcal{D}^{\perp}$  the
result  follows  from  Proposition \ref{prop.n.forms}.
\end{proof}

\medskip\noindent
The following proposition  describes  the equation
$\mathcal{E}_{\mathcal{D}}$ in terms of local coordinates.
\begin{proposition}
Let  $\mathcal{D}\subset\mathcal{C}$ be an n-dimensional
distribution.  Then there  exists a local contact coordinates
$(x^i,z,p_i)$  such that
\begin{equation}\label{eq.local.expr.D}
\mathcal{D}=\langle X_1,X_2,\dots,X_n\rangle\,,\quad   X_i=\widehat{\partial}_{x^i} + b_{ij}\partial_{p_j}
\end{equation}
for some functions $b_{ij}\in C^\infty(M)$. In term of these
coordinates
\begin{equation}\label{eq.coordinates.ED}
\mathcal{E}_{\mathcal D}=\left\{L_P = \langle\widehat{\partial}_{x^i} +
p_{ij}\partial_{p_j}\rangle\,\,\,\big|\,\,\, \det\|p_{ij}-b_{ij}\|=0\,\right\}.
\end{equation}
\end{proposition}
\begin{proof}
The distribution $\mathcal{D}$ can be written in the form
(\ref{eq.local.expr.D}) if
\begin{equation}\label{eq.cond.z}
\mathcal{D}\cap\langle\partial_{p_{1}},\dots,\partial_{p_{n}}\rangle=0.
\end{equation}
Starting  from a local contact system of coordinates $(\bar
x^i,\bar z,\bar p_i)$,  we can construct a new  contact  system of
coordinates of the form
\[
\left\{
\begin{array}
[c]{l}%
{x^{i}}=\bar x^{i}+\epsilon_{i}\bar p_{i}\\
z=\bar z-\frac{1}{2}\sum\epsilon_{i}\bar p_{i}^{2}\\
p_{i}=\bar p_{i}%
\end{array}
\right.
\]
where $\epsilon_i$ are appropriate constants, which satisfies condition \eqref{eq.cond.z}. In terms of these coordinates, the condition
$$
L_P \cap \mathcal{D} =\langle\,\widehat{\partial}_{x^i} +
 p_{ij}\partial_{p_j}\,\rangle\cap \langle\, \widehat{\partial}_{x^i} +
 b_{ij}\partial_{p_j}\rangle \neq 0
$$
is expressed by \eqref{eq.coordinates.ED} in view of Proposition
\ref{prop.local.ED.point}.
\end{proof}

\begin{remark}
The $\omega$-orthogonal complement $\mathcal{D}^\bot$
of $\mathcal{D}$ defines the same equation as $\mathcal{D}:$
$\mathcal{E}_{\mathcal{D}}=\mathcal{E}_{\mathcal{D}^\bot}$.
In general,  the distributions  $\mathcal{D}$ and
$\mathcal{D}^\bot$ are not contactomorphic. As an example, let us
consider the case $n=2$ and the distribution
\[
\mathcal{D}=\langle\widehat{\partial}_{x^1}+x^1\partial_{p_2}\,, \,\, \widehat{\partial}_{x^2}+x^2\partial_{p_1} \rangle.
\]
Its derived distribution
\[
\mathcal{D}'=\langle\widehat{\partial}_{x^1}+x^1\partial_{p_2}\,, \,\, \widehat{\partial}_{x^2}+x^2\partial_{p_1}\,,\,\, \partial_z\rangle
\]
is integrable, whereas the derived distribution of ${\mathcal{D}^{\perp}}$
\[
{\mathcal{D}^{\bot}}'=\langle\widehat{\partial}_{x^1}+x^2\partial_{p_2}\,, \,\, \widehat{\partial}_{x^2}+x^1\partial_{p_1}\,,\,\, (x^2-x^1)\partial_z+\partial_{p_1}-\partial_{p_2}\rangle
\]
is not. In fact it is straightforward to check that $\dim{\mathcal{D}^{\bot}}''=4$.
\end{remark}
In the following theorem, taking into account identification \eqref{eq.rank.1.vectors}, we reformulate the results of Theorem \ref{theorem.main.2}.
\begin{theorem}\label{theorem.main.2.MAE}
Let $m^1\in ({\mathcal{E}_{\mathcal{D}}})_m$ be a
regular point. Then the conformal metric $g_{\mathcal{E}_{\mathcal{D}}}$ is decomposable:
$(g_{\mathcal{E}_{\mathcal{D}}})_{m^1}=\ell_{m^1}\vee\ell'_{m^1}$, where
$\ell_{m^1}=L_{m^1}\cap\mathcal{D}_m$ and
$\ell_{m^1}'=L_{m^1}\cap\mathcal{D}^\bot_m$ are lines. Then there exist only two $(n-2)$-parametric families of characteristic hyperplanes of $L_{m^1}$: one rotates around $\ell_{m^1}$, the other around $\ell_{m^1}'$. Moreover, the characteristic cone is given by
$$
Ch_{m^1}(\mathcal{E}_{\mathcal{D}}) = \{ \pm \eta \otimes \eta,\,\,  \eta \in
\ell^0_{m^1}\cup \ell'^{\,0}_{m^1}\}
$$
where $\ell^0_{m^1},\,\ell'^{\,0}_{m^1} \subset L_{m^1}^* $ are, respectively, the annihilators of $\ell_{m^1}$ and $ \ell'_{m^1}$. Covectors $\eta
\in L_{m^1}^*$ which correspond to characteristic directions and
belong to $\ell^0_{m^1}$ (resp., $\ell'^{\,0}_{m^1}$) define
hyperplanes $\{\eta =0\}$ which contain $\ell_{m^1}$ (resp., $\ell'_{m^1}$).
If one varies the point $m^1$ on
${\mathcal{E}_{\mathcal{D}}}_m$, the line
$\ell_{m^1}$  (resp., $\ell'_{m^1}$) fills the $n$-dimensional
space $\mathcal{D}_{m}$ (resp. $\mathcal{D}_m^\bot$).
\end{theorem}
%
%
Conversely, let us consider a partial differential
equation $\mathcal E\subset M^{(1)}$ which has the following property:
there exists a subdistribution $\mathcal D$ such that for each
$m^1\in\mathcal E$ (over the point $m\in M$),
$$L_{m^1}\cap \mathcal D_m\ne 0.$$
Obviously, in this situation we have that $\mathcal E\subseteq\mathcal
E_{\mathcal D}$. Being both $\mathcal E$ and $\mathcal
E_{\mathcal D}$ submanifolds of the same dimension, locally, they coincide: given $m^1\in\mathcal E$, there exists an
open set $\mathcal{O}\subset M^{(1)}$ containing $m^1$ such that
$$\mathcal E\cap \mathcal{O}=\mathcal E_{\mathcal D}\cap \mathcal{O}.$$

\noindent
This property, without the addition of any other, has no practical value in view of the impossibility of finding the subdistribution
$\mathcal D$. So, in order to have a converse of Theorem \ref{theorem.main.2.MAE}, we have to follow the steps outlined in that theorem.
\begin{theorem}\label{th.criterion}
Let $\mathcal{E} \subset M^{(1)}$ be a $2^{nd}$ order PDE which satisfies  the following properties:
\begin{enumerate}
\item Its conformal metric is decomposable:
$$
(g_{\mathcal{E}})_{m^1} = \ell_{m^1} \vee \ell'_{m^1} $$
where $\ell_{m^1},\ell'_{m^1} \subset L_{m^1} $ are lines.
\item If we let vary the point $m^1$ along the fibre $\mathcal{E}_{m}$, the lines  $\ell_{m^1}, \ell'_{m^1}$ fill two $n$-dimensional spaces ${\mathcal{D}_1}_m$, ${\mathcal{D}_2}_m$ of $\mathcal{C}_{m}$.
\end{enumerate}
Then, locally, $\mathcal{E}=\mathcal{E}_{\mathcal{D}_1}=\mathcal{E}_{\mathcal{D}_2}$.
\end{theorem}

\smallskip\noindent
In the case $n=2$, the above theorem characterizes the classical hyperbolic and
parabolic Monge-Amp{\`e}re equations (i.e. with $2$ independent variables). More precisely we have the following
\begin{corollary}\label{cor.case.n.2}
A second order partial differential equation $\mathcal{E}\subset M^{(1)}$ with $2$ independent variables is a non-elliptic
MAE if and only if the characteristic lines
fill two 2-dimensional subdistributions $\mathcal{D}_1$, $\mathcal{D}_2$ of the contact distribution of $M$. Subdistibutions $\mathcal{D}_1$, $\mathcal{D}_2$ are mutually orthogonal. Moreover, the
equation is parabolic if  $\mathcal{D}_1=\mathcal{D}_1^\perp$ and
is hyperbolic otherwise.
\end{corollary}
\begin{proof} It is sufficient to take into account that, in the case
$n=2$, a MAE $\mathcal{E}$ has characteristic
directions if and only if it is of the form $\mathcal{E}_{\Omega}$
where $2$-form $\Omega$ is decomposable.
\end{proof}

\begin{example}
Let us consider the case $n=2$ and the hyperbolic MAE
\begin{equation}\label{eq.MA.det.ell.n.2}
\mathcal{E}:\,p_{11}p_{22}-p_{12}^2+1=0.
\end{equation}
Equation of characteristics \eqref{eq.characteristics}, restricted to $\mathcal{E}$, is
\begin{equation*}
(p_{12}^2-1)\eta_1^2-2p_{11}p_{12}\eta_1\eta_2+p_{11}^2\eta_2^2=0.
\end{equation*}
The left side term is decomposable in
\begin{equation*}
\big((p_{12}+1)\eta_1-p_{11}\eta_2\big)\big((p_{12}-1)\eta_1-p_{11}\eta_2\big)
\end{equation*}
so that the conformal metric of $\mathcal{E}$ at a point $m^1$ is equal to $(g_{\mathcal{E}})_{m^1}=\ell_{m^1}\vee\ell'_{m^1}$ where
\begin{equation}\label{eq.2.lines.example}
\ell_{m^1}=\langle(p_{12}+1)w_1-p_{11}w_2\rangle\,, \quad \ell'_{m^1}=\langle(p_{12}-1)w_1-p_{11}w_2\rangle
\end{equation}
with
$$
w_1=\widehat{\partial}_{x^1}+p_{11}\partial_{p_1}+p_{12}\partial_{p_2}\,,\quad w_2=\widehat{\partial}_{x^2}+p_{12}\partial_{p_1}+\frac{p_{12}^2-1}{p_{11}}\partial_{p_2}.
$$
Lines \eqref{eq.2.lines.example}  are the only characteristic
subspaces for $\mathcal{E}$ at $m^1$. By a direct computation we
realize that such lines are, respectively
$$
\langle(p_{12}+1)(\widehat{\partial}_{x^1}+\partial_{p_2}) + p_{11}(\partial_{p_1}-\widehat{\partial}_{x^2})\rangle\,, \quad \langle(p_{12}-1)(\widehat{\partial}_{x^1}-\partial_{p_2}) - p_{11}(\partial_{p_1}-\widehat{\partial}_{x^2})\rangle.
$$
If we let vary the point $m^1$ on the fibre $\mathcal{E}_m$, $m=\pi(m^1)$, previous lines fill the following mutually orthogonal $2$-dimensional planes at $m$
$$
\mathcal{D}_m=\langle\widehat{\partial}_{x^1} + \partial_{p_2}\,, \,\,  \widehat{\partial}_{x^2} - \partial_{p_1} \rangle\,, \quad \mathcal{D}_m^\perp=\langle\widehat{\partial}_{x^1} - \partial_{p_2}\,, \,\,  \widehat{\partial}_{x^2} + \partial_{p_1} \rangle
$$
so that we obtain distributions $\mathcal{D}$ and $\mathcal{D}^\perp$ on $M$.

\smallskip\noindent
If we consider two generators of distribution $\mathcal{D}$, for instance $\widehat{\partial}_{x^1} + \partial_{p_2}$ and $\widehat{\partial}_{x^2} - \partial_{p_1}$, we have that
\[
(\widehat{\partial}_{x^1} + \partial_{p_2})\cdot\theta \wedge (\widehat{\partial}_{x^2} - \partial_{p_1})\cdot\theta = dp_1\wedge dp_2 + dp_1\wedge dx^1 + dp_2\wedge dx^2 + dx^1\wedge dx^2
\]
whose restriction on Lagrangian planes gives the $2$-form (see also Remark \ref{rem.horizontalization})
\[
\Omega= (p_{11}p_{22}-p_{12}^2+1)dx^1\wedge dx^2
\]
which vanishes iff Equation \eqref{eq.MA.det.ell.n.2} is
satisfied. We obtain the same result if we consider two generators
of the  distribution $\mathcal{D}^\perp$.
\end{example}

\begin{example}\label{ex.r12}
Let us consider the case $n=3$ and the equation
\begin{equation}\label{eq.MA.example.n.3}
\mathcal{E}:\,\,p_{12}-f(x^i,z,p_i)=0.
\end{equation}
The equation of characteristics \eqref{eq.characteristics} of
$\mathcal{E}$ is $\eta_1\eta_2=0$. Then the conformal metric of $\mathcal{E}$ at a point $m^1$ is equal to $(g_{\mathcal{E}})_{m^1}=\ell_{m^1}\vee\ell'_{m^1}$ where
\[
\ell_{m^1}=\langle\widehat{\partial}_{x^1} + p_{11}\partial_{p_1} + f\partial_{p_2} + p_{13}\partial_{p_3}\rangle\,, \quad \ell_{m^1}'=\langle\widehat{\partial}_{x^2} + f\partial_{p_1} + p_{22}\partial_{p_2} + p_{23}\partial_{p_3}\rangle
\]
If we let vary the point $m^1$ on the fibre $\mathcal{E}_m$, $m=\pi(m^1)$, lines $\ell_{m^1}$ and $\ell'_{m^1}$ fill, respectively, the following mutually orthogonal $3$-dimensional planes at $m$
\[
\mathcal{D}_m=\langle\widehat{\partial}_{x^1} + f\partial_{p_2}, \, \partial_{p_1}, \, \partial_{p_3}\rangle\,, \quad \mathcal{D}^\perp_m=\langle\widehat{\partial}_{x^2} + f\partial_{p_1}, \, \partial_{p_2}, \, \partial_{p_3}\rangle
\]
so that we obtain distributions $\mathcal{D}$ and $\mathcal{D}^\perp$ on $M$.

\smallskip\noindent
If we consider three generators of distribution $\mathcal{D}$, for instance $\widehat{\partial}_{x^1} + f\partial_{p_2}$, $\partial_{p_1}$ and $\partial_{p_3}$, we have that
$$
(\widehat{\partial}_{x^1} + f\partial_{p_2})\cdot\theta\wedge \partial_{p_1}\cdot\theta \wedge \partial_{p_3}\cdot\theta = dp_1\wedge dx^1\wedge dx^3 + fdx^1\wedge dx^2\wedge dx^3
$$
whose restriction on Lagrangian planes gives the $3$-forms (see also Remark \ref{rem.horizontalization})
$$
 \Omega = (-p_{12}+f)dx^1\wedge dx^2 \wedge dx^3
$$
which vanishes iff Equation \eqref{eq.MA.example.n.3} is satisfied. We obtain the same result if we consider three generators of distribution $\mathcal{D}^\perp$.
\end{example}

\section{The full prolongation of a $2^{nd}$ order PDE and its formal integrability}\label{sec.full.prolong}

For the sake of completeness, in this section we consider some
formal aspects of the integration of a $2^{nd}$ order PDE
$\mathcal{E}$. We will treat this subject in the framework of contact manifolds by using, in
addition, the conformal metric $g_{\mathcal{E}}$.

\subsection{The full prolongation  of a contact manifold}
We can define the {\bf $k$-prolongation} $M^{(k)}$ of a contact
manifold $(M, \mathcal{C})$ iteratively as follows. To start with,
we put $M^{(0)}=M$, $\mathcal{C}^{(0)}=\mathcal{C}$ and
$\pi_{1,0}=\pi$. Then we define
\[
M^{(k+1)}=\{\, \text{Lagrangian planes of $M^{(k)}$}  \,\}
\]
where  Lagrangian planes  of $M^{(k)}$ are defined iteratively in
the following way. The manifold $M^{(k)}$ is endowed
with the distribution
\begin{equation}\label{eq.Ck}
\mathcal{C}^{(k)}=\{v\in T_{m^k}M^{(k)}\,\,|\,\,
\pi_{k,k-1\,*}(v)\in L_{m^k}\}
\end{equation}
where $L_{m^k}\equiv m^k$ is a point of $M^{(k)}$ considered as
a Lagrangian plane in $\mathcal{C}^{(k-1)}_{m^{k-1}}$ and
\begin{equation}\label{eq.affine.bundle}
\pi_{k,k-1}:M^{(k)} \to M^{(k-1)}, \,\,\,  m^k \mapsto m^{k-1}
\end{equation}
is the natural projection. It is known \cite{KLR} that \eqref{eq.affine.bundle} are affine bundle for any $k>1$.
Denote by $\theta^{(k)}$  the distribution of $1$-forms on $M^{(k)}$ which defines distribution \eqref{eq.Ck}: $\mathcal{C}^{(k)} = \Ker\,{\theta^{(k)}}$.
\begin{definition}
An $n$-dimensional  subspace $L \subset T_{m^k}M^{(k)}$ is called a {\bf
Lagrangian plane} if it is horizontal w.r.t. $\pi_{k,k-1}$ (i.e. $\pi_{k,k-1\,*}|_{L}$ is not degenerate) and the distributions $\theta^{(k)}$ and $d\theta^{(k)}$ vanish on it.
\end{definition}
In the same way as in Section \ref{subsec.PDE}, a contact chart
$(x^i,z,p_i)$ of $M$  defines  a chart $(x^i, z, p_i, p_{i_1i_2},
\dots,p_{i_1\cdots i_{k+1}})$ of $M^{(k)}$ in a way that a point  $m^k \equiv L_{m^k} \in M^{(k)}$ is given by
$$
L_{m^{k}}=\langle \partial_{x^i} +\sum_{|I| \leq k}
p_{I,i}\partial_{p_{I}} \rangle \,
$$
where $I=(i_1\cdots i_{\ell}),\, 1 \leq i_1\leq i_2 \leq\dots \leq
i_{\ell} \leq n$ is a multi-index of length  $|I|=\ell$   and
$ I,i\overset{\textrm{def}}{=} (i_1,\dots,i_{\ell},i) $ (which
will be reordered if necessary). The distribution $\theta^{(k)}$ is spanned by the 1-forms
$$
\theta_I = dp_I - p_{I,i}dx^i,\quad \vert I\vert\le k.
$$

\smallskip\noindent
Integral manifolds of $\mathcal{C}^{(k)}$ project onto integral manifolds of $\mathcal{C}^{(k-1)}$ through $\pi_{k,k-1}$. In particular, Lagrangian submanifolds $S \subset M^{(k)}$ (i.e. submanifolds such that  $T_sS \in M^{(k+1)}, \forall s \in S$) project onto Lagrangian submanifolds of  $M^{(k-1)}$.

\subsection{The full prolongation of a second order PDE
$\mathcal{E}\subset M^{(1)}$ and its formal integrability}

\smallskip
\noindent The { \bf $1^{st}$-prolongation} of a submanifold
$S\subset M^{(k)}$ is  the submanifold $S^{(1)}\subset M^{(k+1)}$
defined as follows:
\begin{equation*}
S^{(1)}= \text{The set of points $m^{k+1}\in M^{(k+1)}$ such that
$L_{m^{k+1}}$}
\begin{cases}
\subseteq T_{m^k}S\cap\mathcal{C}^{(k)}_{m^k}\quad\text{if $\dim S\ge n$}\\
\supseteq T_{m^k}S\cap\mathcal{C}^{(k)}_{m^k}\quad\text{if $\dim S\le n$}
\end{cases}
\end{equation*}
where $m^k=\pi_{k+1,k}(m^{k+1})$. Iteratively, we define the $h$-prolongation $S^{(h)}\subset M^{(k+h)}$
of $S$.

\smallskip\noindent
We define the {\bf full prolongation} $M^{(\infty)}$ as the
inverse limit of the tower of projections $\dots\longrightarrow
M^{(k)}\overset{\pi_{k,k-1}}{\longrightarrow}
M^{(k-1)}\longrightarrow\dots$ so that a point $m^\infty\in
M^{(\infty)}$ is a sequence $(m=m^0,m^1,\dots,m^k,\dots)$ where $m^k\in M^{(k)}$ and $\pi_{k,k-1}(m^k)=m^{k-1}$. Similarly, we define the full prolongation $S^{(\infty)}$ of any submanifold $S\subset
M^{(k)}$.

\smallskip\noindent
A system of (resp. scalar) PDEs of order $k$, with one unknown
function, is a submanifold (resp. hypersurface) $\mathcal{E}$ of
$M^{(k-1)}$.
\begin{definition}
A \textbf{formal solution} of a $k$-th order PDE $\mathcal{E}$ is a point of $\mathcal{E}^{(\infty)}$.
\end{definition}
Now  we describe  the $k$-th prolongation $\mathcal{E}^{(k)}
\subset M^{(k+1)}$ of a second order PDE
$$
\mathcal{E}=\{F(x^i,z,p_i,p_{ij})=0\}\subset M^{(1)}.
$$
We denote by
$$
D_i = \partial_{x^i} + p_i \partial_z + p_{ij}\partial_{p_j} + \cdots
$$
the total derivative w.r.t. $x^i$  and  for $I=(i_1, \cdots ,i_{\ell})$ we put $D_I = D_{i_1} \circ \cdots \circ D_{i_{\ell}}$. It is straightforward to check that
the $k$-th prolongation
$\mathcal{E}^{(k)}$ of $\mathcal{E}$ is locally described by the system of equations
$$
\mathcal{E}^{(k)} = \{ F=0, \, D_I F =0,\,\,\, 1\leq |I| \leq k   \}.
$$
As a corollary, we can describe the fibre $\mathcal{E}^{(k)}_{m^k} =
\pi_{k+1, k}^{-1}(m^k) \cap \mathcal{E}^{(k)}$ of the projection
$$
\pi_{k,k-1}|_{\mathcal{E}^{(k)}} : \mathcal{E}^{(k)} \to
\mathcal{E}^{(k-1)}
$$
in terms of the  coordinates $p_I$, $|I|=k+2$, of the  fibre $M^{(k+1)}_{m^k} = \pi_{k+1,k}^{-1}(m^k)$ and of the metric
$$
g^{ij}_{\mathcal{E}} =\frac{1}{2-\delta_{ij}}
\frac{\partial F}{\partial p_{ij}}.
$$
We will consider coordinates $p_I = p_{i_1 \cdots i_{\ell}}$ as
symmetric tensor of $S^\ell(\mathbb{R}^n)$.

\begin{corollary} Let $m^1 = (x^i,z,p_i,p_{ij}) \in \mathcal{E}$.
Then $\mathcal{E}^{(1)}_{m^1} $ is defined  by the following
system of linear equations
$$
\mathcal{E}^{(1)}_{m^1} =\{ (2-\delta^{j\ell}) g^{j\ell}p_{ij\ell} = c_i \}
$$
where
$c_i = c_i(m^1)= -\left(\frac{\partial F}{\partial x^i} + p_i \frac{\partial F}{\partial z} +p_{ij}\frac{\partial F}{\partial p_j}\right)(m^1)$. More generally,  if $m^k \in \mathcal{E}^{(k-1)}$, then
$$
\mathcal{E}^{(k)}_{m^k} =\{ (2-\delta^{j\ell})  g^{j\ell}p_{i_1 \cdots i_{k-1}j\ell} = c_{i_1 \cdots i_{k-1}} \}
$$
where
$$ c_{i_1 \cdots i_{k-1}}=c_{i_1 \cdots i_{k-1}}(m^{k-1})=
[D_{i_{k-1}}c_{i_1 \cdots i_{k-2}} -
(D_{i_{k-1}}g^{j\ell})p_{i_1\cdots i_{k-2} j\ell}](m^{k-1}).
$$
\end{corollary}
Recall the following
\begin{definition}
An equation $\mathcal{E} \subset M^{(1)}$ is called
{\bf formally integrable} if the prolongations $\mathcal{E}^{(k)}$ are smooth submanifolds of $M^{(k+1)}$ and
$
\pi_{k+1,k}|_{\mathcal{E}^{(k)}}: \mathcal{E}^{(k)} \to \mathcal{E}^{(k-1)}
$
are smooth fibre bundles.
\end{definition}

\begin{theorem}\label{th.f.i}
Let $\mathcal{E}=\{F=0\} \subset M^{(1)}$ be a smooth hypersurface of $M^{(1)}$. The equation $\mathcal{E}$ is  formally integrable if the associated conformal metric $g_{\mathcal{E}}$ does not vanish (i.e. for any $m^1 \in \mathcal{E}$, $(g_{dF})_{m^1} \neq 0$).
\end{theorem}
To prove the theorem we need the following lemma.
\begin{lemma}
Let  $b =b^{j\ell}(y) \in S^2V^*$ (resp., $c=c_{i_1 \cdots  i_{k-1}}(y) \in S^{k-1}V $) be a symmetric bilinear form (resp., symmetric contravariant $(k-1)$-tensor) in the vector space $V =\mathbb{R}^n  = \{v = (v_1, \cdots v_n)\}$ which smoothly depends on coordinates $y =(y_1, \cdots , y_q) \in \mathbb{R}^q$. If $b\neq 0$ for all $y \in \mathbb{R}^q$, then the equation
\begin{equation}\label{eqcontraction}
b^{j\ell}(y)p_{i_1 \cdots i_{k-1} j \ell} = c_{i_1 \cdots i_{k-1}}(y)
\end{equation}
defines a smooth submanifold $H \subset \mathbb{R}^q\times S^{k+1}V$ such that  the natural projection $\pi : H \to \mathbb{R}^q$ is  an affine fibration with a fibre of  dimension $d(k,n):= \dim S^{k+1}\mathbb{R}^n - \dim S^{k-1}\mathbb{R}^n$.
\end{lemma}
\begin{proof}
First of all, one can easily check  that  the
contraction
$$
\iota_b :  S^{k+1}V \to S^{k-1}V, \,\,p_{i_1 \cdots i_{k-1} j
\ell}\mapsto b^{j\ell}p_{i_1 \cdots i_{k-1} j\ell}
$$
is surjective if $b\neq 0$. This shows that  $\pi^{-1}(y)$ is an affine space of dimension $d(k,n)$. To construct a local
coordinates in $H$, we  consider a linear  change of coordinates
$v_i \to v'_i = A^j_i(y)v_j $ with the matrix $A(y)$ depending on $y$ which transforms the bilinear form $b$ into the  standard form:
$$
b= \epsilon_i \delta^{ij},  \epsilon_i \in \{\pm 1, 0 \}.
$$
We can assume that $\epsilon_1 =1$. The components
$
p_{i_1 \cdots i_{k-1} j\ell}, \, c_{i_1 \cdots i_{k-1}}
$
transform like  tensors. In terms of  the new components
$
p'_{i_1 \cdots i_{k-1} j\ell}, \, c'_{i_1 \cdots i_{k-1}}
$
the equation \eqref{eqcontraction} takes the form
$$
p_{11I} = c_I - \sum_{j>1}\epsilon_j p_{jjI}.
$$
This is a system of linear equations with free variables $p_J,
p_{1J}$ where  the multi-index $J$ does not contain $1$. These
free  variables together with $y$ form a coordinate system of $H$ such that the projection $\pi: H \to \mathbb{R}^q$ is given by
$
\pi(y, p_J, p_{1J})=y.
$
\end{proof}

\smallskip\noindent
\begin{proof}[Proof of Theorem \ref{th.f.i}]
Now we can prove  the theorem by induction. We will assume that
$\mathcal{E}^{(k-1)}\subset M^{(k)}$ is a smooth  submanifold. Then
the restriction of the affine bundle $M^{(k+1)} \to M^{(k)}$ to
$\mathcal{E}^{(k-1)}$ is a locally trivial bundle which  locally
can be identified with the trivial bundle
$$
\mathcal{E}^{(k-1)} \times S^{k+1}\mathbb{R}^n \to \mathcal{E}^{(k-1)},\,  (y, p_{i_1 \cdots i_{k+1}})\mapsto y $$
where $y$ are local coordinates of $\mathcal{E}^{(k-1)}$. Then
$\mathcal{E}^{(k)}$ is defined  by  the  system of equations
$$
(2-\delta^{j\ell}) g^{j\ell}(y) p_{i_1 \cdots i_{k-1} j \ell}  = c_{i_1 \cdots
i_{k-1}}(y)$$
where $g^{j\ell}(y), c_I(y)$ are smooth functions of
$y$. Now the theorem follows from lemma.
\end{proof}

\subsection{Formal solution of a non-characteristic Cauchy problem}\label{sec.formal.solutions}

In this subsection an explicit formal solution of $\mathcal{E}$ is
given once we fix a (non characteristic) Cauchy datum $N$. The
reader can guess that the proof of the following theorem is
related to the possibility of writing the equation $\mathcal{E}$
in the Cauchy-Kowalewski normal form. In fact, this is a
particular instance of a classical result (see for instance \cite{Petro});
a general statement, showing that the existence of
non-characteristic covectors allows to write a system of PDEs in the Cauchy-Kowalewski normal form, was proved in \cite{MMR}.
\begin{theorem}\label{th.form.sol.cauchy}
\bigskip Let $N\subset M$ be an $(n-1)$-dimensional integral manifold of $\mathcal{C}$, $m=m^{0}\in N$, $m^{1}\in\mathcal{E}_{m}$ such that
\begin{equation}
T_{m^{1}}(T_mN)^{(1)}\nsubseteqq T_{m^{1}}\mathcal{E}_{m}.\label{non-char}
\end{equation}
Then, there exists exactly one point $m^{\infty}=\{m^{k}\}_{k\in\mathbb{N}%
_{0}}\in\mathcal{E}^{(\infty)}$ such that, for any $k\in\mathbb{N}_0$, it holds
\begin{equation}
L_{m^{k+1}}\supset T_{m^{k}}N_{\mathcal{E}}^{(k)},\label{Condition}%
\end{equation}
with manifolds $N_{\mathcal{E}}^{(k)}\subset M^{(k)}$ recursively defined by
formulas%
\[
N_{\mathcal{E}}^{(k)}:=(N_{\mathcal{E}}^{(k-1)})^{(1)}\cap\mathcal{E}%
^{(k-1)},~~N_{\mathcal{E}}^{(0)}:=N.
\]
\end{theorem}
Without entering into details, the proof consists in fixing in the
neighborhood of $m$ a Darboux chart $(x^i,z,p_i)$ such that $N$ is represented by%
\begin{equation}\label{eq.fab.1}
\left\{
\begin{array}
[c]{l}%
x^{n}=z=0\\
p_{h}=0,~h<n\\
p_{n}=\Phi_{n}(\widetilde{x})
\end{array}
\right.  ,
\end{equation}
for some suitable function $\Phi_{n}(\widetilde{x})$, $\widetilde{x}=(x_{1},\dots,x_{n-1})$ (see Corollary \ref{cor.caratterizzazione.cauchy.datum}), and showing by a recursive scheme that, in such a chart,
$N_{\mathcal{E}}^{(k-1)}$ is described by
\begin{equation}
\left\{
\begin{array}{l}
x_{n}=z=0 \\
p_I=~\left\{
\begin{array}{cl}
0~ & \text{if \ }i_{a}\leq n-1\;\forall\, a \\
\\
\frac{\partial ^{|J|}}{\partial x^{J}}\Phi_{\underset{h}{\underbrace{n\cdots n}}}(\widetilde{x}) & \text{if \ }I =(J,{\underset{h}{\underbrace{n\cdots n}}}),\,h<\ell\,\,,\,j_b\leq n-1\,\forall\,b \\
\\
\Phi_{\underset{\ell}{\underbrace{n\cdots n}}}(\widetilde{x}) & \text{if \ }I =({\underset{\ell}{\underbrace{n\cdots n}}})
\end{array}
\right.
\end{array}
\right. ,  \label{representation}
\end{equation}
with $\ell$ running from $1$ to $k$, where $I=(i_1\cdots i_\ell)$, $J=(j_1\cdots j_{\ell-h})$, $\partial x^{J}=\partial x^{j_1}\cdots \partial x^{j_{\ell-h}}$ and function $\Phi_{\underset{\ell}{\underbrace{n\cdots n}}}(\widetilde{x})$
is obtained by expliciting jet variable $p_{\underset{\ell}{\underbrace{n\cdots n}}}$ in the equation
$$
\left. (D_{\underset{\ell-2}{\underbrace{n\cdots n}}}F)\right\vert _{(N_{\mathcal{E}}^{(\ell-2)})^{(1)}}=0,
$$
where $\mathcal{E}=\{F=0\}$. This can be done at any step, since the coefficient of the higher order term of $D_{\underset{\ell-2}{\underbrace{n\cdots n}}}F$ (i.e. the coefficient of $p_{\underset{\ell}{\underbrace{n\cdots n}}}$), is $\frac{\partial F}{\partial p_{nn}}(m^1)$, and $\frac{\partial F}{\partial p_{nn}}(m^1)\neq 0$ in view of non-characteristicity
condition \eqref{non-char}. Indeed, let $U=T_{m}N$. By computing the Jacobian matrix of
\eqref{eq.fab.1} one gets
$
U=\langle \xi_{1},\dots,\xi_{n-1}\rangle,
$
with
\begin{equation}
\xi_{h}=\left.\partial_{x^{h}}\right\vert_{m}+\frac{\partial \Phi_{n}}{\partial x^{h}}\left.\partial_{p_{n}}\right\vert_m
 =\left.  \widehat{\partial}_{x^{h}}\right\vert _{m}+\sum_{j=1}^{n}%
p_{hj}(m)\left.  \partial_{p_{j}}\right\vert _{m}
\end{equation}
for $h=1,\dots.n-1$ (with functions $p_{hj}$ given by
\eqref{representation}). But vectors $\xi_{h}$ are exactly the first $n-1$ vectors of the canonical basis of
Lagrangian plane $L_{m^{1}}$, for any $m^{1}\in\pi^{-1}(m)\cap
N^{(1)}$; hence,
$
U^{(1)}=\pi^{-1}(m)\cap N^{(1)}
$
and this curve is described by the free parameter $p_{nn}$, so
that
$
T_{m^{1}}U^{(1)}=\langle\left.
\partial_{p_{nn}}\right\vert _{m^{1}}\rangle;
$
therefore, non-characteristicity condition \eqref{non-char} is
exactly $\frac{\partial F}{\partial p_{nn}}(m^1)\neq 0$.

\smallskip\noindent
Once (\ref{representation}) is proved, it can be
used to check (\ref{Condition}) by simple computations.

\smallskip\noindent
Note that Theorem \ref{th.form.sol.cauchy} is, substantially, an infinitesimal formal analogue of
Cauchy-Kowalewski theorem, and that $m^{\infty}$ corresponds to the Taylor expansion of the unique formal solution of Cauchy problem $(\mathcal{E},N,m)$.

\section{Intermediate integrals of general $2^{nd}$ order PDEs, general MAEs, MAEs of Goursat type and generalized Monge method}\label{sec.int.int.PDEs.MAEs}

\subsection{Intermediate integrals of $2^{nd}$ order PDEs and general MAEs}\label{sec.int.int.PDE}

For the sake of simplicity, we give  the definition of
intermediate integrals only for  PDEs of second order. Recall that
$ M_{f}=\{m\in M\,\,|\,\, f(m)=0\}$  denotes  the zero level set
of  a function  $f\in C^\infty(M)$.

\begin{definition}
Let $\mathcal{E}\subset M^{(1)}$ be a $2^{nd}$  order PDE. A
function $f\in C^\infty(M)$ is called an \textbf{intermediate
integral} of $\mathcal{E}$ if all solutions of $1$-parametric
family $\{M_{f-c}\}_{c\in\mathbb{R}}$ of first order PDEs, are
also solutions of $\mathcal{E}$.
\end{definition}
The following lemma   follows from the  definition of solution of a first order PDE.
\begin{lemma}\label{lem.facile.facile}
A Lagrangian submanifold $\Sigma$ of $M$ is a solution of the first order PDE $f=0$ iff $\Sigma^{(1)}\subset M_f^{(1)}$.
\end{lemma}
We need also the following lemma.
\begin{lemma}\label{lem.facile.facile.2}
Any Lagrangian plane $ L\subset T_mM_f$ is tangent  to a solution
of PDE $M_f$.
\end{lemma}
\begin{proof}
In view of Theorem \ref{th.munoz}, we can suppose that $f=p_n$. Then
$$
 T_mM_f=\langle
\partial_{x^1},\dots,\partial_{x^n},\partial_z,
\partial_{p_1},\dots,\partial_{p_{n-1}}\rangle
$$
and
$$
L=\langle \widehat{\partial}_{x^i} + p_{ij}\partial_{p_j}
\rangle\,,\,\,\, p_{ij}\in\mathbb{R}\,,\,\,p_{nj}=p_{jn}=0.
$$
Now the function
$$
z=z(m)+\sum_{i=1}^{n-1}p_i(m)(x^i-x^i(m)) + \sum_{i,j=1}^{n-1}
p_{ij}\frac{(x^i-x^i(m))(x^j-x^j(m))}{ 2 - \delta_{ij}}
$$
is a solution  tangent to $L$.
\end{proof}

\begin{proposition}
A function f is an intermediate integral of $\mathcal{E}$ iff\,\,
 $\bigcup_{c\in\mathbb{R}} M_{f-c}^{(1)}\subset \mathcal{E}$.
\end{proposition}
\begin{proof}
The condition is necessary. Assume that $f$ is an intermediate
integral.  Let $m^1\equiv L_{m^1}\in M_{f-c}^{(1)}$ for some
$c\in\mathbb{R}$.
 Then by Lemma \ref{lem.facile.facile.2} $m^1$ is tangent to a solution $\Sigma$ of PDE
 $f=c$ which is  also a solution of $\mathcal{E}$. This means  that
$m^1\in\Sigma^{(1)} \subset \mathcal{E}$.

\smallskip\noindent
The condition is sufficient. Let us suppose that $\bigcup_{c\in\mathbb{R}} M_{f-c}^{(1)}\subset \mathcal{E}$. If we fix $c\in\mathbb{R}$,
by Lemma \ref{lem.facile.facile} $\Sigma\subset M$ is solution of
the first order PDE $f=c$ iff $\Sigma^{(1)}\subset
M_{f-c}^{(1)}$,
 which implies that $\Sigma^{(1)}\subset\mathcal{E}$. Hence $\Sigma$ is also a solution of $\mathcal{E}$.
\end{proof}

\begin{theorem}\label{th.f.int.int.Yf.str.char}
A function $f\in C^\infty(M)$ is an intermediate integral of $\mathcal{E}$ iff integral curves of $Y_f$ are strongly characteristic for $\mathcal{E}$.
\end{theorem}
\begin{proof}
Recall that $Y_f=Y_{df}=Y_{f-c}$. Also, $\langle (Y_f)_m \rangle^\perp=\mathcal{C}_m\cap T_mM_{f-f(m)}$. Then ${(Y_f)}_m^{(1)}=({T_mM}_{f-f(m)})^{(1)}$ and theorem follows in view of the above proposition.
\end{proof}

\smallskip\noindent
As an application of previous results we are able to characterize $2^{nd}$ order PDEs which have a large number of intermediate integrals. Such PDEs are described in the following theorem whose statement was known by Goursat \cite{Goursat2}. We give a simple and clear geometric proof of it.
\begin{theorem}\label{th.Goursat.1}
Let $\mathcal{E}$ be a $2^{nd}$ order PDE. If there exist $n$
independent functions $f_1,\dots,f_n$ such that
$f=\varphi(f_1,\dots,f_n)$ is an intermediate integral for any
$\varphi$, then $\mathcal{E}=\mathcal{E}_{\mathcal{D}}$ where
$\mathcal{D}=\langle Y_{f_1},\dots,Y_{f_n} \rangle$.
\end{theorem}
\begin{proof}
For each $f=\varphi(f_1,\dots f_n)$ we have that ${Y}_f^{(1)}\subset \mathcal{E}$ by Theorem \ref{th.f.int.int.Yf.str.char}. Now let us define
$$
\mathcal{D}_m=\{ (Y_f)_m\,\,\,|\,\,\,f=\varphi(f_1,\dots,f_n) \text{  with $\varphi$ arbitrary} \};
$$
it describes an $n$-dimensional subdistribution of $\mathcal{C}$. Indeed, if $\dim\mathcal{D}=n-1$, then $\{Y_{f_1},\dots,Y_{f_n}\}$ would be dependent, and this would imply that the contact form $\theta$ is dependent on $\{df_1,\dots,df_n\}$, which is not possible, as $\theta$ must depend at least on $(n+1)$ differential $1$-forms (see Section \ref{sec.ContactManifolds}).  By definition,
$
\bigcup_{f=\varphi}{(Y_f)}_m^{(1)} = \mathcal{E}_{\mathcal{D}_m}.
$
Since $\bigcup_{f=\varphi}{(Y_f)}_m^{(1)}\subseteq\mathcal{E}_m$, we conclude that $\mathcal{E}_{\mathcal{D}_m}\subseteq\mathcal{E}_m$.
\end{proof}

\smallskip\noindent
The following theorem describes intermediate integrals for any
Monge-Amp{\`e}re equation.
\begin{theorem}[\cite{Alonso Blanco}]\label{th.Ricardo}
A function $f$ is an intermediate integral of a Monge-Amp{\`e}re
equation $\mathcal E_\Omega$, with $\Omega$ an arbitrary $n$-form on
the contact manifold $(M,\mathcal{C})$,
if and only if the associated Hamiltonian vector field $Y_f$
satisfies the following equation:
\begin{equation*}
df\wedge\theta\wedge i_{Y_f}\Omega=0.
\end{equation*}
where $\theta$ is a contact form.
\end{theorem}

\subsection{Intermediate integrals of MAEs of type $\mathcal{E}_{\mathcal D}$}\label{sec.interm.integrals}

Now we describe intermediate integrals for Monge-Amp{\`e}re
equations of type $\mathcal E_\mathcal D$.
\begin{theorem}\label{th.equiv.int.int.with.first.int}
A function $f\in C^\infty(M)$ is an
intermediate integral of the Monge-Amp{\`e}re equation $\mathcal
E_\mathcal D$ if and only if the associated Hamiltonian field
$Y_f$ belongs to $\mathcal D$ or $\mathcal D^\perp$. Equivalently,
the intermediate integrals are the first integrals of
$\mathcal D$ or $\mathcal D^\perp$.
\end{theorem}
\begin{proof}
According to Theorem \ref{th.f.int.int.Yf.str.char}, $f$ is an intermediate integral of
$\mathcal{E}_{\mathcal{D}}$ iff $Y_f$ is strongly characteristic. By arguing as at the beginning of the proof of Theorem \ref{th.D.equal.D.orth}, we obtain that for equations of type $\mathcal{E}_{\mathcal{D}}$ this means that $Y_f\in \mathcal{D}$ or $Y_f\in\mathcal{D}^\perp$.
\end{proof}
\begin{corollary}
If $\mathcal{D}$ (or $\mathcal{D}^\perp$) admits a first integral, or equivalently its derived flag
$$
\mathcal{D}\subseteq \mathcal{D}'\subseteq\mathcal{D}''\subseteq\dots\subseteq\mathcal{D}^k\subseteq\dots
$$
is such that $\mathcal{D}^k\varsubsetneq TM$ for any $k$,
then $\mathcal{E}_{\mathcal{D}}$ admits a smooth solution.
\end{corollary}

\begin{corollary}
The set of intermediate integrals of $\mathcal{E}_{\mathcal{D}}$ is the union of two subrings $\mathcal{R}_1$ and $\mathcal{R}_2$ of $C^\infty(M)$ which are in involution, in the sense that if $f_i\in \mathcal{R}_i$, $i=1,2$, then $\{f_1,f_2\}:=\omega(Y_{f_1},Y_{f_2})=0$.
\end{corollary}

\smallskip\noindent
The following theorem characterizes the simplest equation of type $\mathcal{E}_{\mathcal D}$. Such characterization was known by Goursat \cite{Goursat2}; here we give a proof by using simple properties of contact manifolds together Theorem \ref{th.equiv.int.int.with.first.int}.
\begin{theorem}
The following conditions are equivalent:
\begin{enumerate}
\item $\mathcal{D}$ is an $n$-dimensional integrable distribution of $\mathcal{C}$;
\item $\mathcal{D}$ is generated by $n$ commuting  Hamiltonian vector fields;
\item $\mathcal{E}_{\mathcal D}$ is contact-equivalent to the equation $\det||p_{ij}||=\det||z_{x^ix^j}||=0$;
\item $\mathcal{E}_{\mathcal D}$ is contact-equivalent to the equation $p_{11}=z_{x^1x^1}=0$;
\item $\mathcal{E}_{\mathcal D}$ admits a ring of intermediate integrals generated by $(n+1)$ independent functions.
\end{enumerate}
\end{theorem}
\begin{proof}

\noindent
$1\Rightarrow 2$. In fact, since $\mathcal{D}$ is integrable, we can find $n+1$ functions $\{f_i\}_{i=0...n}$ such that $\mathcal{D}$ is described by $\mathcal{D}=\{df_0=df_1= \cdots  = df_n=0\}$. Since $\mathcal{D}\subset\mathcal{C}$, then (up to a factor)
\begin{equation*}
\theta = df_0 + \sum_{i=1}^{n} a_i df_i
\end{equation*}
for some $a_1,...,a_n\in C^\infty(M)$. Hence
$x^i=f_i$, $z=f_0$,  $p_i=-a_i$, are contact coordinates on $M$ and $\mathcal{D}$ can be written
as
\[
\mathcal{D}=\{dx^1=0,\,\, dx^2=0,\,\dots ,\, dx^n=0,\,\, dz=0\}=\langle \partial_{p_1}\,,\, \dots\,,\partial_{p_n} \rangle.
\]

\smallskip\noindent
$2\Rightarrow 1$. It is an easy application of Theorem \ref{th.munoz}.

\smallskip\noindent
$1\Leftrightarrow 3$. In fact, we already proved that condition $1$ implies that $\mathcal{D}$ is contact-equivalent to $\langle \partial_{p_1}\,, \dots\,,\partial_{p_n}\rangle$. By using  Legendre transformation \eqref{legendre.1} we realize that $\mathcal{D}$ is also contact-equivalent to $\langle \widehat{\partial}_{x^1},\dots, \widehat{\partial}_{x^n} \rangle$, whose associated $\mathcal{E}_{\mathcal D}$ is $\det||p_{ij}||=0$.

\smallskip\noindent
$1\Leftrightarrow 4$. This equivalence goes as the previous one by using a partial Legendre transformation (see \eqref{legendre.2}) which interchanges only $\partial_{p^1}$ with $\widehat{\partial}_{x^1}$.

\smallskip\noindent
$1\Rightarrow 5$. In fact, $\mathcal{D}$ is integrable iff there exist $(n+1)$ functions $f_i$, $i=0,\dots n$, such that $\mathcal{D}=\{df_0=0,\dots,df_n=0\}$. This implies that $\varphi(f_0,f_1,\dots,f_n)$ is a first integral of $\mathcal{D}$ for any function $\varphi$.

\smallskip\noindent
$5\Rightarrow 1$. Let us suppose that $\varphi(f_0,f_1,\dots,f_n)$ is an intermediate integral of $\mathcal{E}_\mathcal{D}=\mathcal{E}_{\mathcal{D}^\perp}$ for any function $\varphi$. In view of Theorem \ref{th.Goursat.1}, $\mathcal{D}$ or $\mathcal{D}^\perp$ is equal to
$
\langle Y_{f_0},...,Y_{f_n} \rangle.
$
Since $\dim \mathcal{D}=n$, then there exist $(n+1)$ smooth functions $\mu_i$, $i=0\dots n$, such that
$$
0=\sum_{i=0}^n \mu_i Y_{f_i}=Y_{\sum \mu_i df_i}
$$
that implies $\sum \mu_i df_i$ depend on the contact form $\theta$, i.e. for some $n$ smooth functions $a_i$ it holds
$$
\theta=df_0+\sum_{i=1}^n a_i df_i
$$
Hence
$x^i=f_i$, $z=f_0$,  $p_i=-a_i$, are contact coordinates on $M$ and $\mathcal{D}$ or $\mathcal{D}^\perp$ can be written
as
\[
\mathcal{D}=\{dx^1=0,\,\, dx^2=0,\,\dots ,\, dx^n=0,\,\, dz=0\}=\langle \partial_{p_1}\,,\, \dots\,,\partial_{p_n} \rangle
\]
which implies that $\mathcal{D}=\mathcal{D}^\perp$.
\end{proof}

\smallskip\noindent

\subsection{Construction of solutions of MAEs of type
$\mathcal{E}_{\mathcal{D}}$ by the generalized Monge method}

As usual, let $(M,\mathcal{C})$ be a contact manifold, $\theta$ a contact form and
$\mathcal{D}\subset\mathcal{C}$ an $n$-dimensional subdistribution of
$\mathcal{C}$. Below we describe a method to construct solutions of
$\mathcal{E}_{\mathcal{D}}$ by generalizing the Monge method of
characteristics (see \cite{Goursat,Morimoto}).
Recall that a vector field $Y\in\mathcal{D}$ is of type $2$ iff
$$
Y\cdot (Y\cdot\theta)=\lambda\theta+\mu (Y\cdot\theta)
$$
for some function $\lambda$ and $\mu$ on $M$.
\begin{proposition}
\label{prop.Monge.method.1}
Let $N\subset M$ be an $(n-1)$-dimensional (embedded) integral
submanifold of the distribution of $\mathcal{C}$ and
$X\in\mathcal{D}$ a vector
field of type $2$ which is transversal to $N$. Let
$$
\Sigma=\bigcup_t\,\varphi_t(N)\subset M
$$
where $\varphi_t$ is the local flow of $X$.
Then $\Sigma$ is solution of the equation $\mathcal{E}_{\mathcal{D}}$ iff
$$
\omega(T_mN,X_m)=0\,\,\forall\,\,m\in N.
$$
\end{proposition}
\begin{proof}
Let us recall that $\Sigma$ is a solution of $\mathcal{E}_{\mathcal{D}}$
if it satisfies the conditions:
\begin{enumerate}
\item $T_m\Sigma\cap\mathcal{D}_m\neq 0\,,\,\,\, \forall\, m\in\Sigma\,$;
\item $T_m\Sigma\subset\mathcal{C}_m\,,\,\,\, \forall\, m\in\Sigma$.
\end{enumerate}
Condition $1$ is obviously satisfied.

\smallskip\noindent
To check condition $2$ we choose coordinates $(t,y^i)$ on $\Sigma$
such that $(y^i)$ are local coordinates on $N$ and $X=\partial_t$.
Any vector field $Y\in\mathcal{X}(N)$ can be considered as
vector field on $\Sigma$ which does not depend on $t$, hence
commutes with $X$. It is sufficient to check
that the function $f(t,y^i):=\theta_{(t,y^i)}(Y)$ be identically zero.
The first two derivatives of $f$ w.r.t. $t$ are
$$
\dot{f}=(X\cdot\theta)(Y)=\omega(X,Y)\,,\,\,\,
\ddot{f}=(X\cdot (X\cdot\theta))Y =
\lambda\theta(Y)+\mu (X\cdot\theta)(Y) = \lambda\,f + \mu\,\dot{f}.
$$
Then $f$ satisfies second order ODE with the initial conditions
$$
f(0,y^i)=0\,,\,\,\,\,\dot{f}(0,y^i)=\omega(X,Y)|_N=0.
$$
This shows that $f\equiv 0$.
\end{proof}

\begin{proposition}\label{prop.alek.1}
Let $\mathcal{D}\subset\mathcal{C}$ be an $n$-dimensional
subdistribution of $\mathcal{C}$. Then a function $f\in
C^\infty(M)$ is a first integral of distribution
$\mathcal{D}^\perp$ ($\mathcal{D}^\perp\cdot f=0)$ iff the
Hamiltonian vector field $Y_f$ belongs to $\mathcal{D}$.
\end{proposition}
\begin{proof}
Let $ \mathcal{D}^\perp=\langle Y_1,\dots,Y_n \rangle $ and $
\mathcal{D}=\Ker\, \theta \cap \Ker\,( Y_1\cdot \theta) \cap
\cdots \cap\Ker\,( Y_n \cdot \theta) .$
The proposition follows from the identity
$$
Y_i(f)=Y_i\,\lrcorner\,df=
Y_f\,\lrcorner\,(Y_i\,\lrcorner d\theta)=
Y_f\,\lrcorner\,(Y_i\cdot\theta).
$$
\end{proof}

\smallskip\noindent
According to Proposition \ref{prop.alek.1}, any first integral $f$
of the distribution $\mathcal{D}^\perp$ defines a Hamiltonian
vector field $Y_f$ (which is a vector field of type $2$, see
Remark \ref{rem.Yf.char.symm}) included in $\mathcal{D}$. So, in
view of Proposition \ref{prop.Monge.method.1}, the problem of
constructing solutions of $\mathcal{E}_{\mathcal{D}}$ reduces to
constructing of $(n-1)$-dimensional submanifolds $N$ of $M$ such
that
\begin{equation}\label{eq.Y.N}
\omega(T_mN,{Y_f}_m)=0\,\,\,\,\forall\,m\in N.
\end{equation}

\begin{proposition}\label{prop.Monge.method.2}
Let $f$ be an intermediate integral of $\mathcal{E}_{\mathcal{D}}$
and $N$ a Cauchy datum for $M_f=\{f=0\}$ (i.e. an $(n-1)$-dimensional integral submanifold of $\mathcal{C}$ included in $M_f$). Then submanifold
$$
\Sigma=\bigcup_t\varphi_t(N)
$$
where $\varphi_t$ is the local flow of Hamiltonian vector field
$Y_f$, is a solution of $\mathcal{E}_{\mathcal{D}}$. If $N$ is
non-characteristic, then the solution is unique.
\end{proposition}
\begin{proof}
Let $X\in TN$. Then
$$
\omega(Y_f,X)=df(X)=X(f)=0.
$$
Therefore $\Sigma$ is a solution of $\mathcal{E}_{\mathcal{D}}$
since $N$ satisfies condition \eqref{eq.Y.N}. The uniqueness of $\Sigma$ follows since $\Sigma$ is also a solution of first order
PDE $f=0$, as it can be derived from Lemma \ref{lemma.char.flow}.
\end{proof}

\smallskip\noindent
Note that $(n-1)$-dimensional submanifold $N\subset M_f$ which is
integral manifold of $\mathcal{C}$ is the same as $(n-1)$-integral
submanifold of the first order PDE
$f\left( x^i,z,{\partial z}/{\partial x^i}\right) =0$.
A description of such submanifolds is given in Section \ref{sec.int.subman.contact}.
In particular, if an $n$-dimensional submanifold $\Sigma$ is a solution of
previous equation, any hypersurface $N$ of $\Sigma$ satisfies above equation.

\smallskip\noindent
Summarizing above results, we can describe a general version of
Monge method of characteristics as follows:
\begin{enumerate}
\item Find a first integral $f$ of the distribution $\mathcal{D}^\perp$. Such
function exists iff $\mathcal{D}^\perp$ belongs to a proper integrable
subdistribution of $TM$. Then the construction of such a function
reduces to finding a solution of a Frobenius system;
\item Find an $(n-1)$-dimensional integral submanifold $N$ of the first order PDE. We can do it by method explained above;
\item Integrate Hamiltonian vector field $Y_f$ to a
local flow $\varphi_t$. Then the submanifold
$$
\Sigma=\bigcup_t\varphi_t(N)
$$
defined in a tubular neighborhood of $N$ is a solution of
$\mathcal{E}_{\mathcal{D}}$.
\end{enumerate}

\begin{theorem}\label{th.Monge.Morimoto}
Let us suppose that $\mathcal{D}$ (or $\mathcal{D}^\perp$) possesses $n$ independent first integrals.
Then any Cauchy datum $N$ can be extended to a solution of $\mathcal{E}_{\mathcal D}$.
\end{theorem}
\begin{proof}
Let $f_1,\dots,f_n$ be independent first integrals of $\mathcal{D}$ (so that any function of them is an
intermediate integral of $\mathcal{E}_{\mathcal D}$). Let denote by $g_i$ the restriction of $f_i$ to $N$.
 Of course the functions $g_i$ are dependent. So there exists a non trivial functional relation
$$
\psi(g_1,\dots,g_n)=0.
$$
The function $f=\psi(f_1,\dots,f_n)$ turns out to be an intermediate integral which vanishes on $N$ and it also satisfies the hypothesis of Proposition \ref{prop.Monge.method.2}. Then the flow of $Y_f$ extends $N$ to a solution of $\mathcal{E}_{\mathcal D}$.
\end{proof}

\begin{theorem}\label{th.ultimo}
Assume that  $\mathcal{D}^\perp$ possesses $n$ independent first integrals $f_1, \cdots, f_n$. Denote by  $M_\mathcal{I}=\bigcup_{\phi} M_{\phi(f_1,\dots,f_n)}^{(1)}$ where $\phi$ is an arbitrary function of $n$ variables. Then
$$
{M}_{\mathcal{I}}=\mathcal{E}_{\mathcal D}.
$$
\end{theorem}
\begin{proof}
${M}_{\mathcal{I}}\subset\mathcal{E}_{\mathcal D}$. In fact,
$L\in {M}_{\mathcal{I}}$ means that $L=T_m\Sigma$, where $\Sigma$ is
a solution of a first order PDE $M_f$ for some fist integral  of
the form $f = \varphi(f_1, \cdots, f_n) $ (such $\Sigma$ exists,
by Lemma \ref{lem.facile.facile.2}). Since $\Sigma$ is also a
solution of $\mathcal{E}_{\mathcal D}$, then
$L\in\mathcal{E}_{\mathcal{D}}$.

\smallskip\noindent
$M_{\mathcal I}\supset\mathcal{E}_{\mathcal D}$.  Let
$L=L_{m^1}\in\mathcal{E}_{\mathcal D}$. Then $L_{m^1}\cap
\mathcal{D}_{\pi(m^1)}$ contains a vector $({Y_f})_{\pi(m^1)}$ for an appropriate
first integral $f$ of $\mathcal{D}^\perp$. As a consequence,
$L\in M_f^{(1)}$.
\end{proof}

\begin{example}
Let $Q$ be a $k$-dimensional
smooth manifold and consider the contact manifold $M:=J^1(Q\times
Q,\mathbb{R})$. Let us take the map
$$
A\colon M=J^1(Q\times Q,\mathbb{R})\to T^*Q\,,\,\,\,
j^1_{q,\overline q}f\mapsto d_{\overline q}i_{q}^*f\,,
$$
where $i_q\colon Q\to Q\times Q$ is defined as $i_q(q')=(q,q')$
for each $q'\in Q$. For each $m\in M$ we define
$\mathcal{D}_m=\Ker A_{*m}\cap\mathcal{C}_m$. In this way we get
an $n$-dimensional subdistribution of $\mathcal C$ (the orthogonal
complement $\mathcal{D}^\perp$ can be also constructed in an
analogous way). If  $x^i$, $\overline{x}^i$ are coordinates on
$Q\times Q$ and $z$ is the coordinate on $\mathbb{R}$, we get a
contact chart $\{x^i,\overline{x}^i, z,p_i,\overline{p}_i\}$. Now,
the local expressions for the subdistributions defined above are
$$\mathcal{D}=\langle \widehat{\partial}_{x^i},\partial_{p_i}\rangle,
 \quad \mathcal{D}^\perp=\langle
 \widehat{\partial}_{\overline{x}^i},\partial_{\overline{p}_i}\rangle.
$$
The Monge-Amp{\`e}re equation $\mathcal{E}_\mathcal{D}$, which is
associated with $2k$-form $\Omega=dp_1\wedge\cdots\wedge
dp_k\wedge dx^1\wedge\cdots\wedge dx^k$, is described in
coordinates by
\begin{equation*}
\textrm{det}\,\left(\frac{\partial^2z}{\partial\overline
{x}^i\partial x^j}\right)=0.
\end{equation*}
Taking into account Theorem \ref{th.equiv.int.int.with.first.int}
and the local expressions of $\mathcal D$ and $\mathcal{D}^\perp$,
the intermediate integrals of $\mathcal{E}_\mathcal{D}$ are
$\varphi(x^1,\dots, x^k,p_1,\dots,p_k)$ and
$\varphi(\overline{x}^1,\dots,
\overline{x}^k,\overline{p}_1,\dots,\overline{p}_k)$, where
$\varphi$ is an arbitrary function of $2k$ variables.

\smallskip\noindent
Therefore, the generalized Monge method applies to
$\mathcal{E}_\mathcal{D}$ and any Cauchy datum can be extended to
a solution in a unique way. In order to illustrate the method we
will carry out all computations in a simple concrete example. Let
$k=2$ so that the equation reads
$$
\frac{\partial^2 z}{\partial \overline{x}^1\partial x^1}
\frac{\partial^2 z}{\partial \overline{x}^2\partial x^2}
-\frac{\partial^2 z}{\partial \overline{x}^1\partial x^2}
\frac{\partial^2 z}{\partial \overline{x}^2\partial x^1}=0.
$$
Now, we consider a Cauchy datum which, for instance, we can
suppose to be parametrizable by $x^1,x^2,\overline{x}^1$; then, we
can fix $\overline{x}^2, p_2$ and $z$ as arbitrary functions of
$x^1,x^2,\overline{x}^1$ and next we determine the remaining
coordinates by imposing the condition of $N$ being a integral
manifold of
$\mathcal{C}=\{dz-p_1dx^1-p_2dx^2-\overline{p}_1d\overline{x}^1-\overline{p}_2d\overline{x}^2=0\}$.
In order to perform explicit computations, let us take, for
example, the Cauchy datum $N$ given by
\begin{equation*}
N\equiv
\begin{cases}
\overline{x}^2=e^{x^2},\,\, p_1=e^{x^1+\overline{x}^1},\,\,
p_2=-{x}^1e^{x^2},\,\, \overline{p}_1=e^{x^1+\overline{x}^1},\,\,
\overline{p}_2=x^1,\,\, z=e^{x^1+\overline{x}^1}
\end{cases}
\end{equation*}
Next, we need to look for an intermediate integral $f=\varphi(x^1,
x^2,p_1,p_2)$ vanishing on $N$. In view of the parametrization of
$N$ we see that $f:=p_2+x^1e^{x^2}$ holds the requirement. The
Hamiltonian field associated with $f$ is
$$
 Y_f=Y_{p_2}+e^{x^2}Y_{x^1}+x^1e^{x^2}Y_{x^2}=
 \partial_{x^2}+p_2\partial_z-e^{x^2}\partial_{p_1}-x^1e^{x^2}\partial_{p_2},
$$
which is easily integrated having the following $8$ first
integrals:
$$
\lambda_1=p_2+x^1e^{x^2},\,\lambda_2=\overline{x}^1,\,
\lambda_3=\overline{x}^2,\,
\lambda_4=\overline{p}_1,\,\lambda_5=\overline{p}_2,\,
\lambda_6=x^1,\, \lambda_7=p_2+x^1p_1,\,\text{and }$$
$$\lambda_8=z-(p_2+x^1e^{x^2})x^2-x^1e^{x^2}.
$$
According with the Theorem \ref{th.Monge.Morimoto}, the propagation of $N$ along the
integral curves of $Y_f$ gives us the unique solution of
$\mathcal{E}_\mathcal{D}$ we are looking for. To do this, it is
sufficient to find $5$ independent relations among the first
integrals of $Y_f$ which hold on $N$, which can be done by
eliminating $7$ coordinates in the parametrization of $N$ by using
the $\lambda$'s. These relations are:
$$
\lambda_1=0,\, \lambda_7=0,\,
\lambda_4-e^{\lambda_6+\lambda_2}=0,\, \lambda_5-\lambda_6=0
\text{ and } \lambda_8-\lambda_4+\lambda_5\lambda_3=0.
$$
By expressing this relations in terms of the original variables we
get, finally,
$$
z=x^1e^{x^2}+e^{x^1+\overline{x}^1}-x^1\overline{x}^2.
$$
\end{example}

\bigskip\noindent\textbf{Acknowledgement.}
This project has been partially supported by RIGS Programme of ICMS, University of Edinburgh. The second author thanks J. Mu\~{n}oz, A. \'{A}lvarez, S. Jim\'{e}nez and J.
Rodr\'{\i}guez for stimulating discussions and encouragements.

\end{document}